\documentclass[pdflatex,sn-mathphys]{sn-jnl}

\usepackage{amsmath}
\usepackage{amsfonts}
\usepackage{amsthm}
\usepackage{mathtools}
\usepackage{bm}

\usepackage{enumitem}
\setlist[enumerate]{leftmargin=.5in}
\setlist[itemize]{leftmargin=.5in}

\DeclareSymbolFont{largesymbolsA}{U}{txexa}{m}{n}
\DeclareMathSymbol{\varprod}{\mathop}{largesymbolsA}{16}

\jyear{2022}%

\theoremstyle{thmstyleone}%
\theoremstyle{thmstyletwo}%
\newtheorem{remark}{Remark}%

\theoremstyle{thmstylethree}%
\newtheorem{myalgorithm}{Algorithm}

\newcommand{\firstreview}[1]{{\color{black} #1}}

\begin{document}

\title{An adaptive hierarchical ensemble Kalman filter with reduced basis models}


\author*[1]{\fnm{Francesco A. B.} \sur{Silva}}\email{f.a.b.silva@tue.nl}

\author[2]{\fnm{Cecilia} \sur{Pagliantini}}\email{cecilia.pagliantini@unipi.it}

\author[1]{\fnm{Karen} \sur{Veroy}}\email{k.p.veroy@tue.nl}

\affil*[1]{\orgdiv{Department of Mathematics and Computer Science}, \orgname{Eindhoven University of Technology}, \orgaddress{\city{Eindhoven}, \postcode{5600 MB}, \country{The Netherlands}}}

\affil[2]{\orgdiv{Department of Mathematics}, \orgname{University of Pisa}, \orgaddress{\city{Pisa}, \postcode{56127}, \country{Italy}}}

\abstract{
The use of reduced order modeling techniques in combination with ensemble-based methods for estimating the state of systems described by nonlinear partial differential equations has been of great interest in recent years in the data assimilation community. Methods such as the multi-fidelity ensemble Kalman filter (MF-EnKF) and the multi-level ensemble Kalman filter (ML-EnKF) are recognized as state-of-the-art techniques. However, in many cases, the construction of low-fidelity models in an offline stage, before solving the data assimilation problem, prevents them from being both accurate and computationally efficient. In our work, we investigate the use of {\it{adaptive}} reduced basis techniques in which the approximation space is modified online by combining information extracted from a limited number of full order solutions and information extracted from reduced models trained at previous time steps. This allows to simultaneously ensure good accuracy and low cost for the employed models and thus improve the performance of the multi-fidelity and multi-level methods.}


\maketitle

\section{Introduction}\label{Introduction}
\label{sec:introduction}

The problem of real time estimation of the state of a physical system known as the {\it{filtering problem}} is common to many applications, from engineering to earth science. Its solution is a main goal of data assimilation (DA) algorithms, which typically combine three sources of information: a deterministic model describing the system's evolution over time, an initial measure of uncertainty on the unknown `true' state, and some updated but polluted information on the state obtained through noisy local measurements. 

The focus of this work is on the state estimation of distributed systems governed by parabolic partial differential equations (PDEs) via synchronous data assimilation, where data are processed in real time, as they are received from the system. This problem requires a framework to quantify and control the uncertainty in the unknown `true' state. In our work, Bayesian statistics serves as the framework of choice. Among the several methods, the sequential Bayesian filter \cite{ho1964} emerges as the most abstract one, serving as a solid foundation for the development of more practical filters. This iterative technique, described in Section \ref{sec:the_filtering_problem}, involves a sequence of steps aimed at updating the measure of the state uncertainty by merging measurements and model information within a Bayesian statistical framework. Although the resulting nonlinear filter is remarkably powerful, it is simply not feasible for most practical applications, as it proves impossible to analytically forecast and update most probability measures. For this reason, other methods, such as the ensemble Kalman filter (EnKF), originally proposed by G. Evensen in \cite{evensen2003}, are commonly preferred. Specifically, the EnKF has demonstrated exceptional robustness in a broad spectrum of nonlinear applications (see, e.g., \cite{counillon2014, sampson2021, diab2021}), despite being based on linear assumptions. Its straightforward implementation and flexibility made it one of the most widespread tools within the data assimilation community.

Over the years, the EnKF has received several enhancements and has evolved into multiple versions, forming what can now be recognized as a `family' of ensemble Kalman filter methods. All these methods use an ensemble of states to represent the uncertainty on the `true' state  and employ a Kalman-like formula for updating the empirical measure of uncertainty. Notably, numerous efforts have been made to improve the method's reliability and effectivity by, for example, introducing inflation terms \cite{anderson1999, hamill2001} to counteract the collapse of covariance due to the finite size of the ensemble \cite{majda2018}, or by employing localization techniques \cite{herschel2002} to enhance the effective ensemble size in the update formula. However, these solutions have only marginal impact on the computational cost of the method. The necessity of solving a large number of instances of the partial differential equations in real time runs remains. In many cases, this is the primary constraint for the effectiveness of the EnKF. 

To address this problem, various solutions have been proposed involving the integration of less costly surrogate models alongside the primary `high-fidelity' model. A notable example can be found in \cite{pagani2017}, where the authors suggest to replace the full order model in the EnKF with a reduced order model (ROM) \cite{hesthaven2022} and to use Kriging \cite{giraldo2011} to build a surrogate for model errors. However,  the reliance on heuristic strategies to manage the approximation error and the challenge in ensuring a priori good approximation properties for the reduced basis space represent major limitations of this method. Notable improvements on the first issue have come from viewing the EnKF as a Monte Carlo method and subsequently integrating state-of-the-art Monte Carlo methodologies into the filtering problem. The approach proposed in \cite{hoel2016} adapts the multi-level Monte Carlo algorithm  \cite{heinrich2001} into a novel multi-level ensemble Kalman filter (ML-EnKF). Likewise, the theory of multivariate control variates \cite{rubinstein85} has been applied to data assimilation  yielding a multi-fidelity ensemble Kalman filter (MF-EnKF) \cite{popov21}.  The idea of using a hierarchy of surrogate models is only one of the possible approaches to improve the approximation properties of the EnKF; other strategies have been proposed in, e.g., \cite{bach2023, narayan2012}.

Both the ML-EnKF and the MF-EnKF combine artificial observations of states simulated at different levels of accuracy to improve the Kalman update. The inclusion of a large number of low accuracy states ensures reliable statistics for the estimated quantities, while the incorporation of a limited number of pairs of states with different levels of accuracy allows to account for the approximation error in an effective and systematic way. However, the selection of suitable surrogate models remains critical for both methods. In the original formulation of the ML-EnKF for spatio-temporal processes \cite{chernov2021}, the surrogate models are assumed to exhibit polynomial convergence, and only spectral methods and multigrid methods are considered for their construction. In contrast, the first implementation of the MF-EnKF also explores the use of reduced basis methods, potentially allowing higher convergence rates. However, the offline construction of these models raises various concerns, particularly when the `true' state is far from the actual dynamics, such as an attractor. In particular, the offline construction of the reduced basis space may not guarantee a sufficient accuracy for an unknown `true' state trajectory that is significantly deviating from the anticipated dynamics. Moreover, the selected space could become too large and costly, especially in the presence of complex dynamics, such as non-periodic and chaotic behaviors.

In \cite{donoghue22} the authors propose a method for dynamically training surrogate models using the trajectories of the high-fidelity states required by the MF-EnKF. The resulting surrogates are locally nearly optimal in a Kolmogorov $n$-width sense \cite{kolmogorov1936}, and the  approximation error is controlled within a fixed tolerance. However, their approach is `memory-less', i.e., it does not retain information on the already built approximation spaces, and, in addition, the admissible error is unrelated to the state uncertainty, which potentially affects the filter's convergence.

In summary, although recent works have attempted to integrate reduced ordered modeling techniques and different extensions of the EnKF method employing multiple levels of accuracy, several issues remain. In particular, the offline construction of a surrogate model appears impractical due to limited knowledge of the system's initial state and the complexity of dynamics over an extended time horizon. The construction of local-in-time reduced basis seems favorable for addressing these issues; however, the size of the training set required to ensure suitable approximation properties remains too large when relying only on local information. Furthermore, to the best of our knowledge, the existing literature lacks a proper definition of adaptive accuracy that effectively links approximation errors to model variability.

In this work we address the issues above by proposing two alternative algorithms, the adaptive reduced basis multi-level ensemble Kalman filter (aRB-ML-EnKF) and the adaptive reduced basis multi-fidelity ensemble Kalman filter (aRB-MF-EnKF).  The proposed algorithms exploit the advantages of an online retraining while  retaining information on the past surrogate models. This approach involves alternating phases of {\it{model inflation}}, in which  the high-fidelity trajectories are used to expand the approximation space, and of {\it{model deflation}}, in which past low-fidelity simulations are used to determine which subspace directions to discard without losing accuracy.

\firstreview{In our study, we focus on parabolic dissipative problems with sufficiently regular solutions. These hypotheses offer an extensively studied framework for the use of reduced-order modeling techniques and ensure a good theoretical understanding of the convergence of data assimilation algorithms \cite{oljava2018}. Furthermore, we consider initial uncertainty distributions that can be effectively approximated within low-dimensional spaces with defined moments and sufficiently light tails, a necessary condition to ensure the problem's reducibility. Within these distributions, we encode the problem's regularity or, when available, some information on the system's long-term behavior.}

The paper is organized as follows. In Section \ref{sec:the_filtering_problem}, we briefly introduce the filtering problem alongside the classical Bayesian and Kalman filters. Then, in Section \ref{sec:reduced_order_modeling}, we offer a brief overview of the technique for reduced order modeling, focusing specifically on proper orthogonal decomposition (POD) and Galerkin projections. Section \ref{sec:multilevel_enkf} is dedicated to a detailed discussion of the ML-EnKF and MF-EnKF, with Subsection \ref{sec:the_prediction_step_of_the_MLEnKF} specifically focused on the innovative elements introduced with our approach. Finally, in Section \ref{sec:numerical_experiments}, we test the proposed approaches on the quasi-geostrophic equations (QGE) and we compare their performances for different surrogate models. In Section \ref{sec:conclusions} we draw conclusions on the proposed method and on its numerical performance.

\section{The filtering problem}
\label{sec:the_filtering_problem}
%


Let $\Omega \subseteq \mathbb{R}^d$, $d\in \{1,2,3\}$, be a spatial domain and let $\mathcal{V} = \mathcal{V}(\Omega)$ be a Hilbert space of finite dimension $\mathcal{N} \in \mathbb{N}^+$ equipped with the inner product $(\,\cdot\,,\,\cdot\,)_\mathcal{V}$. Let $\mathcal{B}(\mathcal{V})$ denote the Borel $\sigma$-algebra on $\mathcal{V}$. The goal of the filtering problem is to update the probability distribution $\pi(\omega(t))$ on $(\mathcal{V}, \mathcal{B}(\mathcal{V}))$, representing the epistemic uncertainty associated with the unknown state $\omega^\star(t)\in\mathcal{V}$, through the assimilation of experimental data.
The data comprises noisy observations obtained at fixed assimilation times, $t_k \in \mathbb{R}$, with  $k \in \{0, ..., N_T\}$ for some $N_T \in \mathbb{N}^+$. The temporal interval between consecutive assimilation times, $\mathcal{I}_k = [t_k, t_{k+1}]$, is referred to as the assimilation window. We denote the noisy observations as $\mathbf{d}_{k}\in \mathcal{Y}$, where $\mathcal{Y}$ is a vector space of size $N_D\in \mathbb{N}^+$, equipped with the inner product $(\,\cdot\,,\,\cdot\,)_\mathcal{Y}$.
In the numerical experiments, the observation space $\mathcal{Y}$ will coincide with $\mathbb{R}^{N_{D}}$.

The measurement process connecting the unknown state $\omega^\star(t_{k})$ at time $t_k$ with the corresponding observations $\mathbf{d}_{k}$ is a stochastic process involving a deterministic measurement operator $\mathcal{L} : \mathcal{V} \rightarrow \mathcal{Y}$, here assumed linear, and an additive noise term $\boldsymbol{\eta}_k \in \mathcal{Y}$. 
Assuming normally distributed noise with zero mean and constant covariance, the measurement model can be expressed as
\begin{align*}
    \mathbf{d}_k = \mathcal{L} \omega ^\star({t}_{k}) + \boldsymbol{\eta}_{k} \quad \text{with} \quad \boldsymbol{\eta}_{k} \sim \mathcal{N}(0, \boldsymbol{\Sigma}),
\end{align*}
where $\boldsymbol{\Sigma} : \mathcal{Y} \rightarrow \mathcal{Y}$ represents the noise covariance, which we assume is known a priori.


The whole measurement process can be presented as the realization of a random process with a random variable distributed according to the likelihood $\pi(\mathbf{d}_k \vert \omega(t_k))$ defined on $(\mathcal{Y}, \mathcal{B}(\mathcal{Y}))$. This distribution represents the probability of obtaining the measurements $\mathbf{d}_k$ given the state $\omega (t_k)$. Knowledge of $\pi(\mathbf{d}_k \vert \omega(t_k))$ is essential for the Bayesian update of the system's uncertainty.
%

The ability to propagate this uncertainty from one point in time to another is also crucial for addressing the filtering problem. In fact, along with the assimilation of new data, it is essential to evolve the state uncertainty consistently with the model expected to describe the evolution of the `true' state. In our case, we assume that the random state, $\omega(t_k) \in \mathcal{V}$, evolves from $t_k$ to $t>t_k$ via the flow \firstreview{$\Psi: \mathcal{V} \times \mathbb{R} \times \mathbb{R} \rightarrow \mathcal{V}$}, so that $\omega (t) = \firstreview{\Psi (\omega (t_{k}), t_k, t)}$. The forecast of the state over each assimilation window $\mathcal{I}_k$ is therefore obtained as
\begin{align}
\label{eq:high_fidelity_dynamical_model}
\omega(t_{k+1}) = \firstreview{\Psi (\omega(t_k), t_k, t_{k+1})}, \qquad \forall k \in \left[ 0, \ldots, N_T \right]. 
\end{align}
%
%
%
For systems governed by partial differential equations, the forecast model \firstreview{$\Psi(\,\cdot\,, t_k, t_{k+1})$} is obtained by integrating a time-dependent PDE over the assimilation window $\mathcal{I}_k$. The general parabolic partial differential equation considered in this work can be formulated in weak form as follows: given the initial condition $\omega(t_k) \in \mathcal{V}$, determine $\omega(t_{k+1}) \in \mathcal{V}$ such that
\begin{align}
    \label{eq:general_ppde} 
    (\partial_t \omega + \mathcal{F}(\omega), v)_\mathcal{V}=0, \qquad \forall v \in \mathcal{V}, \;t \in \left[ t_{k}, t_{k+1} \right].
\end{align}
Here, $\partial_t$ represents the partial derivative in time, and $\mathcal{F}$ is a general differential operator that also encodes boundary conditions. 

To describe the Bayesian filter we extend the
transformation $\firstreview{\Psi(\,\cdot\,, t_k, t_{k+1})}$ to the set of the probability distributions on the measurable space $(\mathcal{V}, \mathcal{B}(\mathcal{V}))$ by introducing the
pushforward operator $\firstreview{\Psi_\#(\,\cdot\,, t_k, t_{k+1})}$.
The Bayesian filter \cite{sarkka2013}, which represents a general probabilistic approach to the filtering problem, can be written in terms of this operation and of Bayesian updates, as in Algorithm \ref{alg:bayes_filter_algorithm}.\\

\begin{myalgorithm} \textbf{The Bayesian filter}
\label{alg:bayes_filter_algorithm}

\noindent
{Input}: a prior probability distribution, $\pi(\omega(t_0))$, characterizing the uncertainty on $\omega^\star({t}_0)$, and the experimental measurements $\{\mathbf{d}_{k}\}_{k=0}^{N_T}$ at times $\{t_k\}_{k=0}^{N_T}$.

\hfill \break
For $k = 0, 1, \ldots, N_T$:\\
\begin{enumerate}[leftmargin=*, label=(\roman*)]
\item
\textbf{Analysis step}:
update the probability distribution according to Bayes' rule
\begin{align*}
    \pi(\omega(t_k) \vert \mathbf{d}_{k}, \mathbf{d}_{k-1}, ...) \propto \pi(\mathbf{d}_k \vert \omega(t_k)) \, \pi(\omega(t_k) \vert \mathbf{d}_{k-1}, \mathbf{d}_{k-2}, ...).
\end{align*}
\item
\textbf{Prediction step}:
predict the probability distribution at time $t_{k+1}$
\begin{align*}
    \pi(\omega(t_{k+1})\vert\mathbf{d}_{k}, \mathbf{d}_{k-1}, ...) = \firstreview{\Psi_\#(\pi(\omega(t_{k})\vert\mathbf{d}_{k}, \mathbf{d}_{k-1}, ...), t_k, t_{k+1})}.
\end{align*}
\end{enumerate}
\end{myalgorithm}

In practice, this filter is impossible to implement in most cases and it must be substituted by approximate filters that may be more practical but possibly less effective. In this work, we focus on the EnKF and on its variants. This family of filters relies on two fundamental assumptions. First, we assume that the probability distributions characterizing the uncertainty on $\omega^\star (t)$ can be approximated by an empirical distribution involving a discrete set of states. Second, we further assume that in the analysis step, all the distributions can be treated as Gaussians and thus, the Kalman update formula can be applied.
The steps of the EnKF and its variants are summarized in Algorithm \ref{alg:enkm_algorithm} and then detailed in Subsection \ref{sec:initial_ensemble} and \ref{sec:the_enkf}.

\begin{myalgorithm} \textbf{The family of ensemble Kalman filters}
\label{alg:enkm_algorithm}

\noindent
{Input}: a prior probability distribution, $\pi(\omega(t_0))$, characterizing the uncertainty on $\omega^\star({t}_0)$, and the experimental measurements $\{\mathbf{d}_{k}\}_{k=0}^{N_T}$ at times $\{t_k\}_{k=0}^{N_T}$.

\hfill \break
(0) \textbf{Initialization}: Generate the initial ensemble 
by sampling states from $\pi(\omega(t_0))$.

\hfill \break
For $k = 0, 1,\ldots, N_T$:\\
\begin{enumerate}[leftmargin=*, label=(\roman*)]
\item
\label{alg:enkm_algorithm:analysis_step}
\textbf{Analysis step}:
update each state in the ensemble(s) with a correction proportional to the misfit between the perturbed state measurements and the experimental data.\\
\item
\label{alg:enkm_algorithm:prediction_step}
\textbf{Prediction step}:
advance each state in the ensemble(s) from time $t_k$ to time $t_{k+1}$ by means of the full order model or a surrogate model.\\
\end{enumerate}
\end{myalgorithm}

In the following sections, we provide a detailed explanation of the steps of Algorithm \ref{alg:enkm_algorithm} for the standard EnKF. Before the data acquisition, the algorithm starts with the initialization of an initial ensemble as described next.


\subsection{Choosing the initial ensembles}
\label{sec:initial_ensemble}

All data assimilation algorithms discussed in this work rely on one or more ensembles to approximate the `true' state uncertainty.  We define the generic ensemble of size $N_\Diamond \in \mathbb{N}^+$ as
\begin{gather*}
E_{k\vert r}^{\Diamond}:=\lbrace \omega_{k\vert r}^{\Diamond,n} \in \mathcal{V} \;\vert\; \omega_{k\vert r}^{\Diamond,n} \sim \pi(\omega(t_k)\vert \mathbf{d}_{r}) \rbrace_{n=1}^{N_\Diamond},
\end{gather*}
where $0\leq k\leq N_T$ and $-1\leq r\leq k$.
The sampling probability distribution $\pi(\omega(t_k)\vert \mathbf{d}_{r})$ corresponds to the uncertainty measure that we aim to approximate using the ensemble.
%
Throughout, we adopt the following convention for the ensemble members: the superscript $\Diamond$ denotes the ensemble family, and 
the index within the ensemble family, i.e., $"\Diamond,n"$ refers to the $n$th state in the $\Diamond$ ensemble; the subscripts $"k\vert r"$ refer to the approximation time $t_k$ and the data assimilation time $t_r$, i.e., $k\vert r$ refers to an approximation of $\omega^\star{(t_k)}$ using data updated at time $t_{r}$. 
The initial ensemble, which we denote by $E_{0\vert -1}^\Diamond$, reflects a partial knowledge of the `true' state $\omega^\star (t_0)$ before any experimental data is made available.
In this work, we consider two alternative sampling distributions for the ensemble $E_{0\vert -1}^\Diamond$, both arising from the knowledge of the model describing the dynamical system. The first option is based on the sole knowledge of
the regularity of the state, e.g., in smooth dissipative processes; the second option assumes that
the dynamics is fully characterized. For example, if the dynamical system possesses an attractor and the `true' state lies within a basin of attraction, the associated invariant measure can be employed as a sampling distribution. In practice, the most common strategy to obtain random states distributed according to this measure involves simulating the system's dynamics for a sufficiently long time and then randomly selecting states from this long-time trajectory. In the numerical experiments of Section \ref{sec:numerical_experiments}, we show how these two possible choices affect the reconstruction of the state.


To account for the $H^1(\Omega)$ regularity of the `true' state $\omega^\star(t)$, as already proposed in \cite{iglesias2016}, we consider a normal distribution with zero mean and covariance $\Delta^{-1}:\mathcal{V}\rightarrow\mathcal{V}$, where $\Delta$ represents the discrete Laplace operator. 
Such a prior effectively privileges smoother modes over less regular ones.
To efficiently sample
from such a normal distribution, 
we generate a set of $\mathcal{N}$ i.i.d. random vectors $\chi^{{\Diamond,n}} \sim \mathcal{N}(0, I)$, with $n=1,...,N_\Diamond$. We then 
compute
\begin{align*}
\omega_{0\vert -1}^{{\Diamond,n}} = \Delta^{-1/2} \chi^{{\Diamond,n}}=  \sum_{i=1}^{\mathcal{N}} \lambda_i^{-1/2}(\xi_i, \chi^{{\Diamond,n}})_\mathcal{V}\,\xi_i \qquad\forall\, n\in\{1,...,N_\Diamond\},
\end{align*}
with $\{\lambda_i\}_{i=1}^\mathcal{N}$ the eigenvalues of the Laplacian and $\{\xi_i\}_{i=1}^\mathcal{N}$ the corresponding eigenfunctions.
\newpage


\subsection{The ensemble Kalman filter}
\label{sec:the_enkf}


The EnKF algorithm starts with the initialization of an initial ensemble, $E^P_{0\vert -1}$, of size $N_P \in \mathbb{N}^+$, whose states are independently sampled according to the rules described in Subsection \ref{sec:initial_ensemble}. 

At the $k$th iteration, with $k\in\{0,\ldots,N_T\}$, the first step of the EnKF is the analysis or update step (Algorithm \ref{alg:enkm_algorithm}, Step {\ref{alg:enkm_algorithm:analysis_step}}) performed when new experimental data are made available. This \firstreview{step} relies on the Kalman update formula to integrate the noisy data $\mathbf{d}_{k}$ with synthetic observations of the predicted states in the forecast ensemble $E^{P}_{k\vert k-1}$.
By assuming an ensemble of states in $\mathcal{V}$, the update step  can be stated as follows: for all $n \in \{1, ..., N_P\}$,
\begin{align*}
    \omega_{k\vert k}^{P, n} = \omega_{k\vert k-1}^{P, n} + \mathbf{K}_k^\text{E} \left( \mathbf{d}_k^{{P,n}} - \mathcal{L} \omega_{k\vert k-1}^{{P,n}} \right)  \quad \text{with} \quad \mathbf{d}^{P,n}_k \sim \mathcal{N}\left( \mathbf{d}_k, \boldsymbol{\Sigma}  \right),
\end{align*}
where $\mathbf{K}_k^\text{E} : \mathcal{Y} \rightarrow \mathcal{V} $ is the Kalman gain \cite{evensen2003} defined as $\mathbf{K}_k^\text{E}  := \mathbf{Q}_k^\text{E} (\mathbf{P}_k^\text{E} + \boldsymbol{\Sigma})^{-1}$. The correlation operators $\mathbf{Q}_k^\text{E} : \mathcal{Y} \rightarrow \mathcal{V}$ and $\mathbf{P}_k^\text{E} : \mathcal{Y} \rightarrow \mathcal{Y}$ are the best unbiased estimators for the state-measurement covariance and for the measurement-measurement covariance, respectively, and can be computed as
\begin{equation*}
\begin{aligned}
\mathbf{Q}^\text{E}_{k} &= \dfrac{1}{N_{P}-1} \sum_{n=1}^{N_{P}} \omega_{k\vert k-1}^{{P,n}} (\mathcal{L} \omega_{k\vert k-1}^{{P,n}}, \,\cdot\,)_\mathcal{Y} - \overline{\omega}_{k\vert k-1}^{P} ( \overline{\mathcal{L}\omega}_{k\vert k-1}^{P}, \,\cdot\,)_\mathcal{Y},\\
\mathbf{P}^\text{E}_{k} &= \dfrac{1}{N_{P}-1} \sum_{n=1}^{N_{P}} \mathcal{L} \omega_{k\vert k-1}^{{P,n}} (\mathcal{L} \omega_{k\vert k-1}^{{P,n}}, \,\cdot\,)_\mathcal{Y} -  \overline{\mathcal{L}\omega}_{k\vert k-1}^{P} ( \overline{\mathcal{L}\omega}_{k\vert k-1}^{P}, \,\cdot\,)_\mathcal{Y},
\end{aligned}  
\end{equation*}
where $\overline{\mathcal{L}\omega}_{k\vert k-1}^{P} \in \mathcal{Y}$ denotes the mean of the ensemble measurements, and
$\overline{\omega}_{k\vert k-1}^P \in \mathcal{V}$ the ensemble mean.

The analysis is followed by the prediction step (Algorithm \ref{alg:enkm_algorithm}, Step \ref{alg:enkm_algorithm:prediction_step}) in which all the states in the post-analysis ensemble $E^P_{k\vert k}$ are advanced over the assimilation window  $\mathcal{I}_{k}$ using the high-fidelity model $\firstreview{\Psi(\,\cdot\,, t_k, t_{k+1})}$, namely 
\begin{equation}
\label{eq:principal_forecast}
    \omega^{P,n}_{k+1\vert k} = \firstreview{\Psi(\omega^{P,n}_{k\vert k}, t_k, t_{k+1})},\qquad \forall\, n \in \{1, ..., N_P\}.
\end{equation}

Solving equation \eqref{eq:principal_forecast} for all ensemble members has arithmetic complexity $\mathcal{O}(\alpha(\mathcal{N}) N_P)$, where $\alpha(\mathcal{N})$ denotes the cost of advancing a single state \firstreview{from time $t_k$ to time $t_{k+1}$}. This is usually performed by solving the PDE \eqref{eq:general_ppde} $N_P$ times and can become prohibitively expensive whenever $\mathcal{N}$ is large.
Hence, this step can become the primary cost in the filtering process due to its dependence on $\mathcal{N}$ and the requirement of a sufficiently large ensemble size $N_P$ to achieve a good data assimilation accuracy.
In this context, the use of surrogate models and reduced order modeling thus becomes crucial.
\newpage

\section{Reduced order modeling}
\label{sec:reduced_order_modeling}
If the initial states and their resulting trajectories belong with high likelihood to a low-dimensional subspace of $\mathcal{V}$, reduced order modeling allows to substitute the full order model with a lower-dimensional surrogate that preserves optimal, or nearly-optimal, approximation properties while significantly diminishing the computational cost required for its solution.
In this study, we concentrate on a specific family of ROM techniques known as {\it{reduced basis methods}} \cite{prudhomme2002}.
The reduced basis method typically consists of two phases: an {\it{offline}} phase and an {\it{online}} phase. In the offline phase, a low-dimensional approximation of the solution space, called the \emph{reduced space}, is identified, and a surrogate model is derived by projecting the full order model onto the reduced space.
The reduced space is constructed by leveraging partial information on the dynamics
derived from training ``snapshots" of the solution at instances in time.
The optimal dimension for this space is determined using a priori error bounds, which allow for the estimation of the introduced approximation error. Then, in the online phase, the reduced-dimensional model resulting from the projection onto the reduced space is solved at a lower computational cost compared to solving the full model. In this context, direct error evaluations and error surrogates, such as a posteriori error bounds, enable the quantification of the approximation error.

In our context, let us consider the solution set
\begin{equation}
\label{eq:continuous_solution_set}
\Xi := \{ \firstreview{\Psi (\omega (t_k), t_k, t)}, \; \forall t \in \mathcal{I},\, \omega(t_k) \sim \pi (\omega (t_k)\vert\mathbf{d}_r)\},
\end{equation}
with $-1\leq r\leq k-1$, obtained by solving the \firstreview{PDE} \eqref{eq:general_ppde} over the interval $\mathcal{I} = [t_k, t_k+T]$ of length $T \in \mathbb{R}^+$.
The problem is considered reducible if the solution set can be accurately approximated by a low-dimensional linear subspace of $\mathcal{V}$, i.e., if, for a given tolerance $\varepsilon  \in \mathbb{R}^+$, there exists a subspace $\mathcal{V}_\varepsilon \subset \mathcal{V}$ of size $N_\varepsilon \ll \mathcal{N}$ such that 
\begin{equation}
\label{eq:continuous_mean_approximation_error}
J(\mathcal{V}_\varepsilon):=\dfrac{1}{T}{ \int_\mathcal{I} \int_\mathcal{V} \big\| (\text{Id} - \Pi_{\mathcal{V}_\varepsilon}) \firstreview{\Psi (\omega(t_k), t_k, t)} \big\|^2_\mathcal{V} \, \pi(d\omega (t_k) \vert \mathbf{d}_r) \,dt }  \leq \varepsilon,
\end{equation}
where $\Pi_{\mathcal{V}_\varepsilon}: \mathcal{V} \rightarrow \mathcal{V}_\varepsilon$ denotes the $\mathcal{V}$-orthogonal projection onto $\mathcal{V}_\varepsilon$ and $\text{Id}:\mathcal{V}\rightarrow\mathcal{V}$ denotes the identity operator.
%
%
Since one typically does not have complete information about $\omega(t_k)$ and its probability distribution, an empirical approximation of the integral above is needed to select a suitable subspace $\mathcal{V}_\varepsilon$, as it will be discussed in Subsection \ref{sec:pod}.

The reduced basis space is built as the span of $N_\varepsilon$ basis functions $\{ \psi_{i}\}_{i=1}^{N_\varepsilon}$ extracted from the snapshots via either SVD-type algorithms, such as the POD \cite{berkooz1993}, or greedy algorithms, or combinations of the two, such as the POD-Greedy method \cite{grepl2005, haasdonk2008}. The high-fidelity solution $\omega(t)$ is then approximated in $\mathcal{V}_{\varepsilon} = \text{span} \{ \psi_{i}\}_{i=1}^{N_\varepsilon} \subset\mathcal{V}$ as
\begin{align*}
\omega_{\varepsilon}(t) = \textstyle \sum_{i=1}^{N_\varepsilon} \omega_{\varepsilon, i}(t) \,\psi_i\,,\quad\;\; (\omega_{\varepsilon, 1}(t),\ldots,\omega_{\varepsilon, N_\varepsilon}(t)) \in\mathbb{R}^{N_\varepsilon},\qquad\forall\, t\in\mathcal{I}.
\end{align*}
The approximate solution is obtained by solving the reduced model: given the initial condition $\omega_\varepsilon (t_k) \in \mathcal{V}_\varepsilon$, determine  $\omega_\varepsilon (t) \in \mathcal{V}_\varepsilon$, for all $t \in \mathcal{I}$, such that
\begin{align}
\label{eq:reduced_ppde}
\left( \partial_t \omega_\varepsilon + \mathcal{F}_{\varepsilon}(\omega_\varepsilon) , v_\varepsilon \right)_\mathcal{V} = 0, \qquad \forall v_\varepsilon \in \mathcal{V}_\varepsilon,
\end{align}
where the operator $\mathcal{F}_{\varepsilon}$ is obtained by projecting $\mathcal{F}$ from \eqref{eq:general_ppde} onto the reduced space $\mathcal{V}_{\varepsilon}$. 
If we consider the problem over the temporal interval $\mathcal{I}_k= [t_k, t_{k+1}]$, we can introduce the reduced flux $\firstreview{\Psi_\varepsilon(\,\cdot\,, t_k, t_{k+1})} : \mathcal{V}_{\varepsilon} \rightarrow \mathcal{V}_{\varepsilon}$ so that
\begin{align*}
\omega_\varepsilon(t_{k+1}) = \firstreview{\Psi_\varepsilon ( \omega_\varepsilon(t_k), t_k, t_{k+1})}.
\end{align*}

Assuming that the full order problem is reducible over the interval $\mathcal{I}_k$ given the prior $\pi(\omega(t_k)\vert \mathbf{d}_k)$,
the reduced model \eqref{eq:reduced_ppde} can be solved at a cost independent of the size $\mathcal{N}$ of the full order problem under the assumption of linearity and parameter-separability of the operator $\mathcal{F}$.
Addressing general non-linear operators may necessitate the application of hyper-reduction techniques, such as empirical interpolation (EIM) \cite{barrault2004}, linear program empirical quadrature \cite{yano2019}, \firstreview{$\ell^p$-quasi-norm minimization \cite{manucci2022},} or empirical cubature \cite{hernandez2017}.

As briefly indicated at the end of Section \ref{sec:introduction}, even when dealing with well-behaved initial distributions, i.e., initial distributions that can be effectively approximated with low-dimensional spaces, in many problems of interest the reducibility of the solution set is not ensured over extended temporal intervals. This entails that large approximation spaces might be needed to recover the state with sufficient accuracy. In this work we propose to rely only on the reducibility of the solution set over brief time horizons by updating the reduced model over time. This allows us to incorporate information on the updated uncertainty, have approximation spaces of smaller dimension and, hence, save on computational time.

\firstreview{The idea of using online adaptive reduced-order modeling techniques to lower the offline training costs and address problems with limited global reducibility is relevant and of current interest. Notable examples can be found in the fields of optimization \cite{carlberg2008, peherstorfer2015, qian2017}, simulation \cite{washabaugh2012, pagliantini2021}, and control \cite{kramer2017} of dynamical systems, as well as in parameter estimation \cite{kartmann2023} and uncertainty quantification \cite{vidlickova2022}. Unlike existing methods, we neither use non-intrusive approaches nor develop evolution equations for the approximation spaces. Instead, we propose a discontinuous online adaptation of the reduced basis, in line with the filter's discontinuous nature. Also, unlike trust region methods, we do not rely on error bounds for space inflation but use a POD-based approach. Additionally, we introduce a deflation phase to remove unused bases, thus adapting the size of the approximation space to the problem's local reducibility.}

We propose an adaptive strategy which entails that the training phase is no longer confined to an offline phase, as in classical reduced basis methods. Instead, we employ online generated snapshots coming from the data assimilation algorithm itself. This approach substantially reduces the computational costs but also introduces extra challenges, which will be discussed and addressed in Subsection \ref{sec:the_prediction_step_of_the_MLEnKF}.
\newpage


\subsection{Proper orthogonal decomposition}
\label{sec:pod}

Proper orthogonal decomposition has shown significant effectiveness in handling parabolic problems \cite{kunisch2001}, and it has been used to construct reduced basis models on-the-fly in a multi-fidelity data assimilation framework (see \cite{donoghue22}). 
%
In this section, we illustrate the connection between POD and the reducibility condition presented in Equation \eqref{eq:continuous_mean_approximation_error}. Specifically, we describe how a discrete solution set can be leveraged to construct a quasi-optimal reduced basis space, with average accuracy within the tolerance $\varepsilon$, and how the optimization problem can be tracked back to an eigenvalue problem on the covariance of this discrete solution set.
The idea is to construct an approximation space via POD and to perform the Galerkin projection of the PDE onto this space.

Ideally, for a fixed size $N_\varepsilon$, the optimal approximation space (in the sense of the Kolmogorov $n$-width \cite{kolmogorov1936, pinkus2012}) would be obtained by minimizing $J(\mathcal{V}_\varepsilon)$, with $\text{dim} \, \mathcal{V}_\varepsilon = N_\varepsilon$. However, determining this optimal space necessitates knowledge of the solution set $\Xi$, defined in Equation \eqref{eq:continuous_solution_set}, which is unattainable in practice. 
An approximation can be obtained by means of Monte Carlo integration, employing a set of snapshots, called the training set, and defined as 
\begin{gather*}
\Xi_\Delta := \{\firstreview{\Psi (\omega^{n}(t_k), t_k, t_k + p \Delta t)},\;\forall p \in \{0, ..., N_\Theta\}, \,\omega^{n}(t_k) \in E_{k\vert r} \},
\end{gather*}
where $E_{k\vert r} := \{ \omega^{n}(t_k) \sim \pi(\omega(t_k)\vert \mathbf{d}_r) \}_{n=1}^{N_\Xi}$, $N_\Xi \in \mathbb{N}^+$, is an ensemble of random initial states, and $\Delta t := T / N_\Theta$, $N_\Theta \in \mathbb{N}^+$, is a time step, sufficiently small to capture the system's dynamics.
Quasi-optimal approximation spaces can be constructed by minimizing the cost function
\begin{gather*}
J_\Delta({\mathcal{V}}_\varepsilon) := 
\dfrac{1}{N_\Xi (N_\Theta + 1)} \sum_{n=1}^{N_\Xi} \sum_{p=0}^{N_\Theta} \big\| (\text{Id} - \Pi_{{\mathcal{V}}_\varepsilon})\omega^n(t_k+p\Delta t) \big\|_\mathcal{V}^2,
\end{gather*}
with $\omega^n(t_k+p\Delta t) \in \Xi_\Delta$. The solution of the resulting optimization problem is closely related to the eigenvalue problem associated with the \firstreview{linear operator} $\mathcal{R}:\mathcal{V}\rightarrow\mathcal{V}$, defined as
\begin{gather*}
\mathcal{R}:= \dfrac{1}{N_\Xi (N_\Theta+1)} \sum_{n=1}^{N_\Xi} \sum_{p=0}^{N_\Theta} \omega^n(t_k+p\Delta t)(\omega^n(t_k+p\Delta t), \,\cdot\,)_\mathcal{V}.
\end{gather*}
Let $\gamma_1 \geq ... \geq \gamma_\mathcal{N} \geq 0$ denote the eigenvalues of $\mathcal{R}$ with associated eigenfunctions $\psi_1, ..., \psi_\mathcal{N} \in \mathcal{V}$. Then, the first $N_\varepsilon$ eigenfunctions $\{\psi_i \}_{i=1}^{N_\varepsilon}$ form a $\mathcal{V}$-orthogonal basis for the optimal subspace $\mathcal{V}_\varepsilon$, while the sum of the remaining eigenvalues, $\{\gamma_i\}_{i>N_\varepsilon}$, provides a measure of the average approximation error, namely
\begin{align}
\label{eq:error_bound}
\min J_\Delta (\mathcal{V}_\varepsilon) = \sum_{i=N_\varepsilon+1}^{\mathcal{N}} \gamma_{i}.
\end{align}
In practice, the eigenpairs of the operator $\mathcal{R}$ can be determined by leveraging the eigenpairs of the Gramian matrix $\mathbf{G} \in \mathbb{R}^{N_\Xi(N_\Theta+1) \times N_\Xi(N_\Theta+1)}$, with elements $(\omega^n(t_k+p\Delta t), \omega^m(t_k+q\Delta t))_\mathcal{V}$, $p, q \in \{0,...,N_\Theta\}$, $n,m \in \{1, ..., N_\Xi\}$, as proven in \cite[Proposition 1]{kunisch2001}. 
Hereafter, we denote the entire algorithm for constructing the reduced basis space $\mathcal{V}_\varepsilon$ from the training set $\Xi_\Delta$ and the training tolerance $\varepsilon$ as $\mathcal{V}_\varepsilon = \text{POD} (\Xi_\Delta, \varepsilon)$.

Leveraging the result expressed in Equation \eqref{eq:error_bound}, we can efficiently identify the optimal size ensuring an approximation error below a specified threshold, i.e., for a given $\varepsilon_k \in \mathbb{R}^+$, we can select the smallest $N_\varepsilon$ for which $\min J_\Delta(\mathcal{V}_\varepsilon) \leq \varepsilon_k$. However, determining an appropriate value for \firstreview{$\varepsilon_k$} remains an open problem. \firstreview{Ideally, the value should be lower than the average distance between states sampled from $\pi(\omega(t_k)\vert \mathbf{d}_r)$ to avoid confusing approximation error and state uncertainty. Based on this intuition, in the following section, we propose a new tolerance selection strategy that links $\varepsilon_k$ to the covariance of $\pi(\omega(t_k)\vert \mathbf{d}_r)$.} 


Let $\omega(t_k)$ and ${\omega'}(t_k)$, be two statistically independent random variables distributed according to the same probability law $\pi(\omega(t_k)\vert \mathbf{d}_k)$, and let $\varepsilon_r \in (0,1)$ be a user-defined relative tolerance for the projection error, then we define the local tolerance as
\begin{equation}
\begin{aligned}
\label{eq:tolerance}
\varepsilon_k :=& \;\varepsilon_r \int_\mathcal{V} \int_\mathcal{V} \big\| \omega(t_k) - {\omega'}(t_k) \big\|_\mathcal{V}^2 \, \pi(d{\omega'}(t_k)\vert \mathbf{d}_k) \, \pi(d\omega(t_k)\vert \mathbf{d}_k)  \\
=&\; 2 \varepsilon_r \int_\mathcal{V} \big\| \omega(t_k) -  \overline{\omega}(t_k) \big\|_\mathcal{V}^2 \, \pi(d\omega(t_k)\vert \mathbf{d}_k) = 2 \varepsilon_r \, \text{tr} \left( \mathbf{C}^{\omega, \omega}_k \,\cdot\,, \,\cdot\,  \right)_\mathcal{V},
\end{aligned}
\end{equation}
where $\mathbf{C}^{\omega, \omega}_k : \mathcal{V} \rightarrow \mathcal{V}$ denotes the covariance of $\omega(t_k) \sim \pi(\omega(t_k), \mathbf{d}_k)$, and $\overline{\omega}(t_k) \in \mathcal{V}$ its mean.
If the exact probability measure is approximated by an empirical measure, as it happens in the EnKF algorithms, the covariance can be approximated by an appropriate Monte Carlo estimator, as shown in Subsection \ref{sec:the_prediction_step_of_the_MLEnKF}.

This definition implies that for a given measure of uncertainty, the distance between a random state $\omega(t_k) \in \mathcal{V}$ and its projection in $\mathcal{V}_\varepsilon$ is on average $\varepsilon_r$ times smaller than the average distance between two random states. This result is particularly relevant in a data assimilation context as it allows the surrogate model to become increasingly accurate as the state uncertainty decreases.

\section{Adaptive hierarchical ensemble Kalman filter with reduced basis models}
\label{sec:multilevel_enkf}

As mentioned in Section \ref{sec:introduction}, a significant improvement to the Kalman update step and to the approximation properties of the EnKF has been achieved with multi-level and multi-fidelity strategies.
Both the ML-EnKF \cite{popov21} and the MF-EnKF \cite{hoel2016} follow the steps of Algorithm \ref{alg:enkm_algorithm}, but, unlike the classical EnKF,
they employ advanced Monte Carlo methods in the analysis step, enabling the use of states from multiple ensembles with different levels of accuracy. Furthermore, they both require a hierarchy of models to simulate, in the prediction step, the evolution of ensembles at different levels of accuracy. 

In this work, we propose to generate these models via adaptive \firstreview{reduced order modeling} of the full model \eqref{eq:high_fidelity_dynamical_model}. 
This operation is conducted dynamically, employing a limited number of costly expensive training data obtained from the full order model alongside several inexpensive data from surrogate models trained in previous steps. With respect to the use of training data from the full order model alone, this
allows for better surrogate models that can retain information and yield more precise reconstructions.
For the sake of simplicity, we focus on a two-level setting, and we refer to Subsection \ref{sec:multilevel} for a brief description of the generalization to an arbitrary number of levels.

In the two-level setup, three ensembles are employed: the {\it{principal}} ensemble $E^P_{k\vert r}$, the {\it{ancillary}} ensemble $E^A_{k\vert r}$, and the {\it{control}} ensemble $E^C_{k\vert r}$, according to the notation introduced in \cite{popov21}. These ensembles have size $N_P, N_A, N_C \in \mathbb{N}^+$, respectively, with the principal and control ensembles, strongly correlated, of identical size, $N_P = N_C$. At the $k$th assimilation time, the principal ensemble comprises states in $\mathcal{V}$, while the other two ensembles exclusively contain states in $\mathcal{V}_\varepsilon^k$. In contrast, during the initialization, all ensemble states are sampled in $\mathcal{V}$, with the principal and ancillary ensembles being sampled independently, i.e., $E^P_{0\vert -1} \perp E^A_{0\vert -1}$, and the principal and the control ensembles being identical, i.e., $E^P_{0\vert -1} = E^C_{0\vert -1}$. 

In the next sections we describe the analysis step of the ML-EnKF (Subsection \ref{section:analysis_ML-EnKF}) and of the MF-EnKF (Subsection \ref{section:analysis_MF-EnKF}), while in Subsection \ref{sec:the_prediction_step_of_the_MLEnKF} we propose a novel strategy for the prediction step.

\subsection{The ML-EnKF analysis step}
\label{section:analysis_ML-EnKF}
In contrast with the standard EnKF, the analysis step of the ML-EnKF updates the states in the ensembles with the empirical Kalman gain
$\mathbf{K}_k^\text{ML} := \mathbf{Q}^\text{ML}_k (\mathbf{P}^\text{ML}_k + \boldsymbol{\Sigma})^{-1}$ where the state-measurement covariance $\mathbf{Q}^\text{ML}_{k} : \mathcal{Y} \rightarrow \mathcal{V}$ and the measurement-measurement covariance $\mathbf{P}^\text{ML}_{k} : \mathcal{Y} \rightarrow \mathcal{Y}$ are obtained via multi-level Monte Carlo estimation.
This operation involves (1) the telescopic expansion of the two covariances in terms of partial covariances associated with models of different accuracy and ensembles of different sizes, and (2) a regularization step to ensure the stability of the method. Considering only the three ensembles introduced in the previous section, and two levels of accuracy associated with the approximation spaces $\mathcal{V}$ and $\mathcal{V}_\varepsilon^k$, the unregularized covariances $\widetilde{\mathbf{Q}}^\text{ML}_{k}: \mathcal{Y} \rightarrow \mathcal{V}$ and $\widetilde{\mathbf{P}}^\text{ML}_{k}: \mathcal{Y} \rightarrow \mathcal{Y}$ can be defined as
\begin{equation}\label{eq:uncov}
\begin{aligned}
\widetilde{\mathbf{Q}}^\text{ML}_k &:= \mathbf{Q}^{P,P}_{k} - \mathbf{Q}^{C,C}_k + \mathbf{Q}^{A,A}_{k},\qquad\\
\widetilde{\mathbf{P}}^\text{ML}_k &:= \mathbf{P}^{P,P}_{k} - \mathbf{P}^{C,C}_k + \mathbf{P}^{A,A}_{k}.
\end{aligned}
\end{equation}
This definition involves the single ensemble state-measurement covariances $\mathbf{Q}^{P,P}_{k} : \mathcal{Y} \rightarrow \mathcal{V} $, and $\mathbf{Q}^{C,C}_{k}, \mathbf{Q}^{A,A}_{k} : \mathcal{Y} \rightarrow \mathcal{V}_\varepsilon^k \subset \mathcal{V}$, computed as
\begin{equation}
\begin{aligned} 
\label{eq:state-measurement-covariances}
\mathbf{Q}^{\Diamond,\Diamond}_{k} = \dfrac{1}{N_{\Diamond}-1} \sum_{n=1}^{N_{\Diamond}} \omega_{{k}\vert {k-1}}^{{\Diamond,n}} (\mathcal{L}\omega_{{k}\vert {k-1}}^{{\Diamond,n}},\,\cdot\,)_\mathcal{Y} - \overline{\omega}_{k\vert k-1}^{\Diamond} ( \overline{\mathcal{L}\omega}_{k\vert k-1}^{\Diamond}, \,\cdot\,)_\mathcal{Y},\qquad
\end{aligned}
\end{equation}
with $\Diamond\in\{P,C,A\}$, and the single ensemble measurement-measurement covariances
$\mathbf{P}^{\Diamond,\Diamond}_{k} : \mathcal{Y} \rightarrow \mathcal{Y}$ defined as
\begin{equation}
\begin{aligned}
\label{eq:measurement-measurement-covariances}
\mathbf{P}^{\Diamond,\Diamond}_{k} = \dfrac{1}{N_{\Diamond}-1} \sum_{n=1}^{N_{\Diamond}} \mathcal{L}\omega_{{k}\vert {k-1}}^{{\Diamond,n}} ( \mathcal{L}\omega_{{k}\vert {k-1}}^{{\Diamond,n}}, \,\cdot\,)_\mathcal{Y} - \overline{\mathcal{L}\omega}_{k\vert k-1}^{\Diamond} (\overline{\mathcal{L}\omega}_{k\vert k-1}^{\Diamond}, \,\cdot\,)_\mathcal{Y}.
\end{aligned}
\end{equation}
%

Although, by construction, the covariance matrices are symmetric and positive semi-definite, their combination in \eqref{eq:uncov} may
yield an indefinite symmetric covariance $\widetilde{\mathbf{P}}^\text{ML}_k$. The presence of a negative part of the spectrum has the potential to trigger instabilities in the data assimilation process. \firstreview{To counteract this issue, regularization techniques are commonly used, such as projecting the matrices on the cone of positive semi-definite matrices via spectrum thresholding \cite{bickel2008_1, chernov2021, hoel2016}, banding or tapering the matrix \cite{anderson2012, bickel2008, olivier2004}, or performing the combination in a log-Euclidean geometry \cite{maurais23}.}

In this work, we propose a regularization based on the projection of both $\widetilde{\mathbf{P}}^\text{ML}_k$ and $\widetilde{\mathbf{Q}}^\text{ML}_k$ to the cone of positive semi-definite linear operators, so as to move the unstable measurement subspace to their kernel. 
More precisely, let $(\lambda_{i}, \mathbf{p}_{i}) \in \mathbb{R} \times \mathcal{Y}$, for all $i \in \{1, ..., N_D\}$, be the eigenpairs of $\widetilde{\mathbf{P}}^\text{ML}_k$. The regularized covariances $\mathbf{P}^\text{ML}_k : \mathcal{Y} \rightarrow \mathcal{Y}$ and $\mathbf{Q}^\text{ML}_k : \mathcal{Y} \rightarrow \mathcal{V}$ are computed as follows:
\begin{equation*}
\begin{aligned} 
\mathbf{Q}^\text{ML}_k &:= \sum_{\substack{{i}=1 \\ \lambda_{i} \geq 0}}^{N_{{D}}} \widetilde{\mathbf{Q}}^\text{ML}_k \mathbf{p}_{i} ( \mathbf{p}_{i} , \,\cdot\,)_\mathcal{Y},
\;\\
\mathbf{P}^\text{ML}_k &:= \sum_{\substack{{i}=1 \\ \lambda_{i} \geq 0}}^{N_{{D}}} \widetilde{\mathbf{P}}^\text{ML}_k \mathbf{p}_{i} ( \mathbf{p}_{i}, \,\cdot\,)_\mathcal{Y} = \sum_{\substack{{i}=1 \\ \lambda_{i} \geq 0}}^{N_{{D}}} \lambda_{{i}} \mathbf{p}_{i} ( \mathbf{p}_{i}, \,\cdot\,)_\mathcal{Y}.
\end{aligned}
\end{equation*}
We remark that our proposed approach is a variation of the regularization proposed in \cite{chernov2021}, with the difference that we project both the state-measurement and the measurement-measurement covariance. Empirically, this choice has shown better stability properties compared to the original regularization scheme \cite{chernov2021}. We also explored the approach of \cite{hoel2016}, which consists in the regularization of the state-state covariance. Although this method exhibited comparable stability properties, it is impractical due to the need to decompose high-dimensional covariances.

Finally, the regularized covariances are used to define the multi-level Kalman gain $\mathbf{K}_k^\text{ML} := \mathbf{Q}^\text{ML}_k (\mathbf{P}^\text{ML}_k + \boldsymbol{\Sigma})^{-1}$, which is employed to update the states of the three ensembles. Specifically, for all $n \in \{1,...,N_P\}$ and $m \in \{1,...,N_A\}$, we update the ensemble members as
\begin{equation*}
\begin{aligned}
{\omega}_{k\vert k}^{{P,n}} &=  \omega_{k\vert k-1}^{{P,n}} + \;\;\;\; \mathbf{K}^\text{ML}_k \left( \mathbf{d}_k^{{P,n}} - \mathcal{L} \omega_{k\vert k-1}^{{P,n}} \right) & &\text{with} \quad \mathbf{d}^{P,n}_k \sim \mathcal{N}\left( \mathbf{d}_k, \boldsymbol{\Sigma}  \right),\\
{\omega}_{k\vert k}^{{C,n}} &=  \omega_{k\vert k-1}^{{C,n}} + \Pi_{\mathcal{V}_\varepsilon^{k}} \mathbf{K}^\text{ML}_k \left( \mathbf{d}_k^{{C,n}} - \mathcal{L} \omega_{k\vert k-1}^{{C,n}}\right) & &\text{with} \quad \mathbf{d}^{C,n}_k = \mathbf{d}^{P,n}_k,\\
{\omega}_{k\vert k}^{{A,m}} &=  \omega_{k\vert k-1}^{{A,m}} + \Pi_{\mathcal{V}_\varepsilon^{k}} \mathbf{K}^\text{ML}_k \left( \mathbf{d}_k^{{A,m}} - \mathcal{L} \omega_{k\vert k-1}^{{A,m}} \right) & &\text{with} \quad \mathbf{d}^{A,m}_k \sim \mathcal{N}\left( \mathbf{d}_k, \boldsymbol{\Sigma}  \right),
\end{aligned}
\end{equation*}
where $\Pi_{\mathcal{V}_\varepsilon^{k}}:\mathcal{V} \rightarrow \mathcal{V}_\varepsilon^k$ is the orthogonal projection from the high-fidelity to the low-fidelity space. This operation guarantees that the updated control and ancillary states belong to $\mathcal{V}^k_\varepsilon$.

\subsection{The MF-EnKF analysis step}
\label{section:analysis_MF-EnKF}

Similarly to the ML-EnKF, the MF-EnKF computes the Kalman gain based on a combination of the cross-covariance between states and state measures of multiple ensembles with different accuracy and ensemble sizes. Unlike the ML-EnKF, the structure of the two covariances $\mathbf{Q}^\text{MF}_{k}$ and $\mathbf{P}^\text{MF}_{k}$ is based on the theory of multivariate control variates \cite{rubinstein85}, a Monte Carlo approach intended to reduce the variance of the estimates. Considering the same three ensembles as before, $E^{P}_{k\vert r}$, $E^{C}_{k\vert r}$ and $E^{A}_{k\vert r}$, and the same two levels of accuracy, we define the two covariances as
\begin{equation}
\begin{aligned}
\label{eq:multi-fidelity-matrices}
\mathbf{Q}^\text{MF}_k &= \mathbf{Q}^{P,P}_k + \dfrac{1}{4} \left( \mathbf{Q}^{C,C}_k + \mathbf{Q}^{A,A}_k \right) - \dfrac{1}{2} \left( \mathbf{Q}^{P,C}_k + \mathbf{Q}^{C,P}_k \right), \\
\mathbf{P}^\text{MF}_k &= \mathbf{P}^{P,P}_k + \dfrac{1}{4} \left( \mathbf{P}^{C,C}_k + \mathbf{P}^{A,A}_k \right) - \dfrac{1}{2} \left( \mathbf{P}^{P,C}_k + \mathbf{P}^{C,P}_k \right)\firstreview{,}
\end{aligned}
\end{equation}
\firstreview{where}, the state-measurement covariances, $\{\mathbf{Q}^{\Diamond,\Diamond}_{k}\}_{\Diamond\in\{P,C,A\}}$,
and the measurement-measure\-ment covariances, $\{\mathbf{P}^{\Diamond,\Diamond}_{k}\}_{\Diamond\in\{P,C,A\}}$, associated with a single ensemble are defined as in \eqref{eq:state-measurement-covariances} and \eqref{eq:measurement-measurement-covariances}, respectively. The remaining state-measurement cross-covariances $\mathbf{Q}^{C,P}_{k} : \mathcal{Y} \rightarrow \mathcal{V}^k_\varepsilon \subset \mathcal{V}$ and $\mathbf{Q}^{P,C}_{k} : \mathcal{Y} \rightarrow \mathcal{V} $ are computed as
\begin{equation}
\begin{aligned} 
\label{eq:missing_state-measurement-covariances}
\mathbf{Q}^{C,P}_{k} &= \dfrac{1}{N_{P}-1} \sum_{n=1}^{N_{P}} \omega_{{k}\vert {k-1}}^{{C,n}} ( \mathcal{L}\omega_{{k}\vert {k-1}}^{{P,n}} , \,\cdot\,)_\mathcal{Y} - \overline{\omega}_{k\vert k-1}^{C} ( \overline{\mathcal{L}\omega}_{k\vert k-1}^{P}, \,\cdot\,)_\mathcal{Y},\\
\mathbf{Q}^{P,C}_{k} &= \dfrac{1}{N_{P}-1} \sum_{n=1}^{N_{P}}  \omega_{{k}\vert {k-1}}^{{P,n}} ( \mathcal{L}\omega_{{k}\vert {k-1}}^{{C,n}} , \,\cdot\,)_\mathcal{Y} - \overline{\omega}_{k\vert k-1}^{P} ( \overline{\mathcal{L}\omega}_{k\vert k-1}^{C} , \,\cdot\,)_\mathcal{Y},
\end{aligned}
\end{equation}
and the measurement-measurement cross-covariances $\mathbf{P}^{C,P}_{k}, \mathbf{P}^{P,C}_{k} : \mathcal{Y} \rightarrow \mathcal{Y}$ as 
\begin{equation}
\begin{aligned} 
\label{eq:missing_measurement-measurement-covariances}
\mathbf{P}^{C,P}_{k} &= \dfrac{1}{N_{P}-1} \sum_{n=1}^{N_{P}} \mathcal{L}\omega_{{k}\vert {k-1}}^{{C,n}} ( \mathcal{L}\omega_{{k}\vert {k-1}}^{{P,n}} , \,\cdot\,)_\mathcal{Y} - \overline{\mathcal{L}\omega}_{k\vert k-1}^{C} ( \overline{\mathcal{L}\omega}_{k\vert k-1}^{P}, \,\cdot\,)_\mathcal{Y},\\
\mathbf{P}^{P,C}_{k} &= \dfrac{1}{N_{P}-1} \sum_{n=1}^{N_{P}} \mathcal{L}\omega_{{k}\vert {k-1}}^{{P,n}} ( \mathcal{L}\omega_{{k}\vert {k-1}}^{{C,n}} , \,\cdot\,)_\mathcal{Y} - \overline{\mathcal{L}\omega}_{k\vert k-1}^{P} ( \overline{\mathcal{L}\omega}_{k\vert k-1}^{C}, \,\cdot\,)_\mathcal{Y}.
\end{aligned}
\end{equation}
%
The multi-fidelity Kalman gain is then defined as $\mathbf{K}_k^\text{MF} := \mathbf{Q}^\text{MF}_k (\mathbf{P}^\text{MF}_k + \tfrac{1}{2} \boldsymbol{\Sigma})^{-1}$, with $\mathbf{P}^\text{MF}_k$ and $\mathbf{Q}^\text{MF}_k$ as in \eqref{eq:multi-fidelity-matrices}. The states of the three ensembles are updated, for all $n \in \{1,...,N_P\}$ and $m \in \{1,...,N_A\}$, as
\begin{equation*}
\begin{aligned}
\widetilde{\omega}_{k\vert k}^{{P,n}} &=  \omega_{k\vert k-1}^{{P,n}} + \;\;\;\;\;\; \mathbf{K}^\text{MF}_k \left( \mathbf{d}_k^{{P,n}} - \mathcal{L} \omega_{k\vert k-1}^{{P,n}} \right) & &\text{with} \quad \mathbf{d}^{P,n}_k \sim \mathcal{N}\left( \mathbf{d}_k, \dfrac{1}{2} \boldsymbol{\Sigma}  \right),\\
\widetilde{\omega}_{k\vert k}^{{C,n}} &=  \omega_{k\vert k-1}^{{C,n}} + \Pi_{\mathcal{V}_\varepsilon^{k}} \mathbf{K}^\text{MF}_k \left( \mathbf{d}_k^{{C,n}} - \mathcal{L} \omega_{k\vert k-1}^{{C,n}}\right) & &\text{with} \quad \mathbf{d}^{C,n}_k = \mathbf{d}^{P,n}_k,\\
\widetilde{\omega}_{k\vert k}^{{A,m}} &=  \omega_{k\vert k-1}^{{A,m}} + \Pi_{\mathcal{V}_\varepsilon^{k}} \mathbf{K}^\text{MF}_k \left( \mathbf{d}_k^{{A,m}} - \mathcal{L} \omega_{k\vert k-1}^{{A,m}} \right) & &\text{with} \quad \mathbf{d}^{A,m}_k \sim \mathcal{N}\left( \mathbf{d}_k, \dfrac{1}{2} \boldsymbol{\Sigma}  \right).
\end{aligned}
\end{equation*}
Unlike the ML-EnKF, the analysis step of the MF-EnKF requires a further re-centering phase, as detailed in \cite{popov21}. This is meant to ensure that the control and ancillary ensembles share the same mean and that the control ensemble remains strongly correlated with the principal ensemble. For this purpose, we introduce the re-centering mean 
\begin{gather*}
\overline{\omega}_{k\vert k}^{\text{MF}} := \overline{\omega}_{k\vert k}^{P} + \dfrac{1}{2} ( \overline{\omega}_{k\vert k}^{A} - \overline{\omega}_{k\vert k}^{{C}})\;\;\mbox{where}\;\;
\overline{\omega}_{k\vert k}^{\Diamond} = \dfrac{1}{N_{\Diamond}} \textstyle \sum_{n=1}^{N_{\Diamond}} \tilde{\omega}_{k\vert k}^{{\Diamond,n}},\;\;
\Diamond\in\{P,C,A\}.
\end{gather*}
The re-centering mean is substituted to the mean of the three ensembles obtained after the Kalman update imposing the following translations:
\begin{equation*}
\begin{aligned}
\omega_{k\vert k}^{P,n} &= \widetilde{\omega}_{k\vert k}^{P,n} - \overline{\omega}_{k\vert k}^{P} + \overline{\omega}_{k\vert k}^{\text{MF}}, &\qquad \forall\, n \in \{1,...,N_{P}\}, \\
\omega_{k\vert k}^{C,n} &= \widetilde{\omega}_{k\vert k}^{C,n} - \overline{\omega}_{k\vert k}^{C} + \Pi_{\mathcal{V}_\varepsilon^k} \overline{\omega}_{k\vert k}^{\text{MF}}, &\qquad \forall\, n \in \{1,...,N_{P}\}, \\
\omega_{k\vert k}^{A,m} &= \widetilde{\omega}_{k\vert k}^{A,m} - \overline{\omega}_{k\vert k}^{A} + \Pi_{\mathcal{V}_\varepsilon^k} \overline{\omega}_{k\vert k}^{\text{MF}}, &\qquad \forall\, m \in \{1,...,N_{A}\}.
\end{aligned}
\end{equation*}




\begin{remark}
\label{remark_1}
In the mean-field limit, as the number of states in the three ensembles tends to infinity, the mean-field ML-EnKF converges to the mean-field EnKF. This is directly related to the convergence of $\mathbf{Q}_k^\text{ML}$ and $\mathbf{Q}_k^{E}$ to the measurement-measurement mean field covariance $\overline{\mathbf{Q}}_k : \mathcal{Y} \rightarrow \mathcal{V}$, and of $\mathbf{P}_k^\text{ML}$, $\mathbf{P}_k^{E}$ to  the state-measurement mean field covariance $\overline{\mathbf{P}}_k : \mathcal{Y} \rightarrow \mathcal{Y}$. The same result does not hold for the MF-EnKF, whose mean-field covariances deviate from the ones of the EnKF unless $\mathcal{V}=\mathcal{V}_\varepsilon^k$. For this reason, unlike what was shown for the mean-field ML-EnKF in \cite{hoel2016}, it is not possible to relate the MF-EnKF results to the EnKF ones.
\end{remark}

\subsection{The prediction step: an inflation-deflation adaptive reduced basis approach}
\label{sec:the_prediction_step_of_the_MLEnKF}

In this section, we introduce our novel approach that combines adaptive reduced order modeling techniques with the prediction phases of the ML-EnKF and of the MF-EnKF algorithm. The primary objective of this approach is to obtain an accurate estimation of the system's `true' state at minimal computational cost; a feat not currently realized by existing state-of-the-art methods. For this purpose, in the prediction step, we dynamically construct surrogate models using information derived from the trajectories of the principal ensemble and the reduced space inherited from the preceding iteration. 
In contrast to global reduced order modeling approaches, this method allows to exploit the local reducibility of the problem over time. Moreover, the introduction of a memory term enables an increase in the effective training set size, implicitly involving the past trajectories. This results in improved approximation properties and lower-dimensional surrogate models.


At each iteration of the EnKF algorithm, we construct two subspaces of $\mathcal{V}$: the first one $\mathcal{V}_\varepsilon^{k}$, called the \emph{inflated space}, is the reduced basis space; the second one $\mathcal{W}_\varepsilon^k$, called the \emph{deflated space}, is the approximation space that transfers information from one assimilation window to the subsequent one. 
Before starting the data assimilation process, along with the initial ensembles, it is necessary to specify the starting deflated space, $\mathcal{W}_{\varepsilon}^{0}$, to incorporate any potential prior knowledge regarding the model's reducibility. In the absence of more information, one can assume $\mathcal{W}_\varepsilon^{0} \equiv \mathcal{V}$. 
Next, at the $k$th iteration of the data assimilation algorithm, the prediction step starts with the forecast of the post-analysis principal ensemble, $E_{k\vert k}^P$, over the assimilation window $\mathcal{I}_{k}$. Specifically, the high-fidelity model $\firstreview{\Psi(\,\cdot\,, t_k, t_{k+1})}$ is employed to update the states in the principal ensemble
\begin{equation}
    \label{eq:forecast_principal_ensemble}
    \omega_{k+1\vert k}^{P,n} = \firstreview{\Psi(\omega_{k\vert k}^{P,n}, t_k, t_{k+1})} ,\quad \forall\, n \in \{1, ..., N_P\}.
\end{equation}

Assuming that the time-stepping used to compute the model action, i.e., to discretize the parabolic PDE \eqref{eq:general_ppde}, requires $N_S \in \mathbb{N}^+$ sub-steps, the computation of the forecast ensemble in \eqref{eq:forecast_principal_ensemble} generates a set of intermediate high-fidelity states at times $t_k + s (t_{k+1} - t_k)N_s^{-1}$, for $1\leq s\leq N_S$, that we collect in the set
$\mathcal{E}_k^P:= \{ \omega^{P, n}_{k+s N_S^{-1}\vert k},\,1\leq n\leq N_P,\, 1\leq s\leq N_S \} $.
This set is used to expand the deflated space $\mathcal{W}_\varepsilon^{k}$ inherited from the previous iteration of the algorithm.

Let $\Xi_{\Delta,k}^{\text{HF}} = \mathcal{E}_k^P - \Pi_{\mathcal{W}_\varepsilon^{k}} \mathcal{E}_k^P$ represent the projection of $\mathcal{E}_k^P$  onto the orthogonal complement of $\mathcal{W}_\varepsilon^{k}$. This set carries information about the dynamics of the principal ensemble over the time window $\mathcal{I}_{k}$ which is not captured by $\mathcal{W}_\varepsilon^{k}$. Consequently, we inflate the subspace $\mathcal{W}_\varepsilon^{k}$ with the POD basis extracted from the training set $\Xi_{\Delta,k}^\text{HF}$, prescribed the tolerance $\varepsilon^\Diamond_k/2$, i.e., we define
\begin{gather}
\label{eq:inflation_pod}
    \mathcal{V}_\varepsilon^{k+1} := \text{span} \{\text{POD} (\Xi_{\Delta,k}^\text{HF}, \varepsilon^\Diamond_k/2),  \mathcal{W}_\varepsilon^{k} \},
\end{gather}
where $\varepsilon^\Diamond_k$ is computed from the post-analysis ensembles $\{E_{k\vert k}^{\Diamond}\}_{\Diamond\in\{P,C,A\}}$ and is based on the filter of choice.
In particular, in the case of the aRB-ML-EnKF, $\varepsilon_k^\text{ML} \in \mathbb{R}$ is computed using a multi-level estimator of the trace for the state covariance, consistently with Equation \eqref{eq:tolerance}. This results in the absolute tolerance
\begin{equation}
\label{eq:multilevel_tolerance}
\begin{aligned}
    \varepsilon^\text{ML}_k = 
    \dfrac{2 \varepsilon_r}{N_P-1} \sum_{n=1}^{N_P}
    \Big(\| \omega_{k\vert k}^{P,n} - \overline{\omega}_{k\vert k}^{P} \|_\mathcal{V}^2 -
    \| \omega_{k\vert k}^{C,n} - \overline{\omega}_{k\vert k}^{C} \|_\mathcal{V}^2\Big)\\
    +\, \dfrac{2 \varepsilon_r}{N_A-1} \sum_{m=1}^{N_A}  \| \omega_{k\vert k}^{A,m} - \overline{\omega}_{k\vert k}^{A} \|_\mathcal{V}^2,
\end{aligned}
\end{equation}
where $\varepsilon_r$ is a given relative tolerance.
In the case of the aRB-MF-EnKF, a multi-fidelity estimator is used for this purpose, leading to the absolute tolerance $\varepsilon_k^\text{MF} \in \mathbb{R}$ computed as
\begin{equation}
\label{eq:multifidelity_tolerance}
\begin{aligned}
    \varepsilon^\text{MF}_k = \, \dfrac{2\varepsilon_r}{N_P-1} \sum_{n=1}^{N_P}
    \Big(\| \omega_{k\vert k}^{P,n} - \overline{\omega}_{k\vert k}^{P} \|_\mathcal{V}^2
    - (\omega_{k\vert k}^{P,n} - \overline{\omega}_{k\vert k}^{P},  \omega_{k\vert k}^{C,n} - \overline{\omega}_{k\vert k}^{C})_\mathcal{V}\\
    +\dfrac14 \| \omega_{k\vert k}^{C,n} - \overline{\omega}_{k\vert k}^{C} \|_\mathcal{V}^2\Big)
    + \dfrac12 \dfrac{\varepsilon_r}{N_A-1}\sum_{n=1}^{N_A} \| \omega_{k\vert k}^{A,m} - \overline{\omega}_{k\vert k}^{A} \|_\mathcal{V}^2.
\end{aligned}
\end{equation}
%
These correspond to a multi-level Monte Carlo estimation, Equation \eqref{eq:multilevel_tolerance}, and to a control variate method estimation, Equation \eqref{eq:multifidelity_tolerance}, of the trace of the covariance of state uncertainty at time $t_k$, consistently with what is outlined in Equation \eqref{eq:tolerance}.

The resulting approximation space $\mathcal{V}_\varepsilon^{k+1}$ of size $N_\varepsilon^{k+1}$ is used to construct the surrogate model $\firstreview{\Psi^{k+1}_\varepsilon(\,\cdot\,, t_k, t_{k+1})}$, following the procedure of Section \ref{sec:reduced_order_modeling}. The updated model is then employed to advance the states in the control and ancillary ensembles, namely
\begin{equation}
\label{eq:forecast_control_and_ancillary_ensemble}
\begin{aligned}
    \omega_{k+1\vert k}^{C,n} & = \firstreview{\Psi^{k+1}_\varepsilon(\Pi_{\mathcal{V}_\varepsilon^{k+1}}\omega_{k\vert k}^{P,n}, t_k, t_{k+1})}  ,&\quad \forall\,n\in \{1,...,N_P\},\\
    \omega_{k+1\vert k}^{A,m} & = \firstreview{\Psi^{k+1}_\varepsilon( \Pi_{\mathcal{V}_\varepsilon^{k+1}}\omega_{k\vert k}^{A,m}, t_k, t_{k+1})},&\quad\forall\,m\in \{1,...,N_A\}. 
\end{aligned}
\end{equation}
The use of $\Pi_{\mathcal{V}_\varepsilon^{k+1}}\omega_{k\vert k}^{P,n}$ as the initial condition for forecasting the control ensemble is essential to maintain the principal and the control ensemble strongly correlated. This choice makes unnecessary the updating of the control ensemble in the analysis step.

Just like for the principal ensemble, the forecast of the control and ancillary ensembles \eqref{eq:forecast_control_and_ancillary_ensemble} generates two sets of 
low-fidelity states denoted as 
%
\begin{align*}
\mathcal{E}_k^C&:= \{ \omega^{n, C}_{k+sN_S^{-1}\vert k}, 1\leq n\leq N_P,\, 1\leq s\leq N_S\},\;\,\\
\mathcal{E}_k^A&:= \{ \omega^{m, A}_{k+sN_S^{-1}\vert k} , 1\leq m\leq N_A,\, 1\leq s\leq N_S\},
\end{align*}
respectively. These are used to deflate $\mathcal{V}_\varepsilon^{k+1}$, retaining only the modes actively employed by the surrogate model: if $\Xi_\Delta^\text{LF} := \mathcal{E}_k^C \cup \mathcal{E}_k^A$, the updated deflated space is
\begin{align}%
\label{eq:deflation_pod}
    \mathcal{W}_\varepsilon^{k+1} := \text{POD} (\Xi_\Delta^\text{LF}, \varepsilon^\Diamond_k/2),
 \end{align}
where $\varepsilon^\Diamond_k$ is defined according to Equation \eqref{eq:multilevel_tolerance} or \eqref{eq:multifidelity_tolerance} depending on the filter of choice.
This operation is computationally less expensive than model inflation, as it solely involves states belonging to a subspace of size $N_\varepsilon^{k+1} \ll \mathcal{N}$. This observation reflects the complete decoupling of this reduced order modeling operation from the size of the full order problem.

\firstreview{The rationale behind the choice of the POD tolerance in Equation \ref{eq:inflation_pod} and Equation \ref{eq:deflation_pod} is to build an approximation space with a prescribed tolerance $\varepsilon^\Diamond_k$ by means of two intermediate approximation spaces whose error contributions could be orthogonal. Intuitively, by ensuring that each of the two spaces introduces an average error equal to half the tolerance $\varepsilon^\Diamond_k$, we aim at obtaining a final model error that is lower than or equal to $\varepsilon^\Diamond_k$. This is, however, only an empirical choice. In fact, the approximation error associated with the space constructed in the deflation phase refers to a different time window from the one in which the inflated model is designed to operate. This involves an additional error contribution that is difficult to estimate.}

\begin{remark}\label{rmk:yano}
    We point out that the aRB-MF-EnKF has connections with existing literature. Indeed, if we discard the deflated space, i.e., $\mathcal{W}_\varepsilon^{k} = \emptyset$, for all $k \in \{0, ..., N_T\}$, the resulting memory-less method coincides with the approach introduced in \cite{donoghue22}. Moreover, by neglecting the entire retraining process and maintaining the approximation space introduced at the initial time $t_0$, i.e., $\mathcal{V}_\varepsilon^k = \mathcal{W}_\varepsilon^0$, for all $k \in \{0, ..., N_T\}$, the reduced order modeling approach simplifies to the one employed in \cite{popov21}.
\end{remark}

\subsection{Telescopic extension}
\label{sec:multilevel}

In this section, we discuss the extension from the two-level to the multi-level scenario for both the multi-level and the multi-fidelity algorithms. The analytical aspects of this operation for the analysis step have been extensively discussed in  \cite{hoel2016} and \cite{popov21}, to which we refer. Consequently, our emphasis is specifically on the prediction step of the two algorithms.

Let $\ell \in \{0,...,L\}$, $L \in \mathbb{N}^+$ represent a series of accuracy levels associated with both the inflated spaces $\mathcal{V}_\varepsilon^{k, \ell} \subseteq \mathcal{V}$ and the deflated spaces $\mathcal{W}_\varepsilon^{k, \ell} \subseteq \mathcal{V}$. These follow a nested inclusion relationship, $\mathcal{V}_\varepsilon^{k,0} \subseteq ... \subseteq  \mathcal{V}_\varepsilon^{k,L-1} \subseteq\mathcal{V}$, and $\mathcal{W}_\varepsilon^{k,0} \subseteq ... \subseteq \mathcal{W}_\varepsilon^{k,L-1} \subseteq \mathcal{V}_\varepsilon^{k,L-1}$, with $\mathcal{V}_\varepsilon^{k,L}=\mathcal{W}_\varepsilon^{k,L}=\mathcal{V}$, for all $k \in \{0, ..., N_T \}$.
This hierarchy of spaces corresponds to a sequence of models $\firstreview{\Psi_\varepsilon^{k,\ell}(\,\cdot\,, t_k, t_{k+1})} : \mathcal{V}^{k,\ell} \rightarrow \mathcal{V}^{k,\ell}$ and of paired ensembles $(E_{k\vert k}^{P,0}, E_{k\vert k}^{C,0} \equiv \emptyset)$, $(E_{k\vert k}^{P,\ell}, E_{k\vert k}^{C,\ell})$, defined as follows
\begin{equation*}
\begin{aligned}
    E_{k\vert k}^{P,\ell}&:=\{\omega_{k\vert k}^{P,n,\ell} \in \mathcal{V}_\varepsilon^{k,\ell} \}_{n=1}^{M_\ell}, \quad &&\forall \ell \in \{0, ..., L\}\\ 
    E_{k\vert k}^{C,\ell}&:=\{\omega_{k\vert k}^{C,n,\ell} = \Pi_{\mathcal{V}_\varepsilon^{k,\ell-1}} \omega_{k\vert k}^{P,n,\ell} \in \mathcal{V}_\varepsilon^{k,\ell-1}\}_{n=1}^{M_\ell},  \quad &&\forall \ell \in \{1, ..., L\}.
\end{aligned}
\end{equation*}
Here $M_\ell$ denotes the size of the ensembles at the $\ell$th accuracy level, for any $\ell \in \{0,...,L\}$. 
The principal ensemble at the $0$th level coincides with the ancillary ensemble of the two-level case. 
\firstreview{The size of each ensemble should be determined considering the error and costs of the surrogate models. Ideas on how to manage this operation and maximize the efficiency of multi-level and multi-fidelity Monte Carlo methods can be found in \cite{jakeman2020, oljava2018, peherstorfer2016,  schaden2020}.}

Similarly to the two-level case, the ensemble forecast in this multi-level scenario yields discrete solution sets $\mathcal{E}_k^{P,\ell}$ and $\mathcal{E}_k^{C,\ell}$ for all $\ell \in \{0,...,L\}$, with $\mathcal{E}_k^{C,0}=\emptyset$. These sets contain the intermediate states required for the model integration over the time window $\mathcal{I}_k$. They can be used for the model inflation and deflation, i.e., to construct the inflated spaces
\begin{gather*}
    \mathcal{V}_\varepsilon ^{k+1,\ell} = \text{span} \{\text{POD}(\Xi_{\Delta, k}^{\text{HF}}, \varepsilon_k^{\Diamond, \ell} /2 ),  \mathcal{W}_\varepsilon^{k,\ell} \}, \qquad \forall \ell \in\{0,...,L-1\},
\end{gather*}
from the high-fidelity training set $\Xi_{\Delta, k}^{\text{HF}} = \mathcal{E}_k^{P, L} -  \Pi_{\mathcal{W}_\varepsilon^{k,\ell}} \mathcal{E}_k^{P,L}$, and the deflated spaces
\begin{gather*}
    \mathcal{W}_\varepsilon ^{k+1,\ell} = \text{POD}(\Xi_{\Delta,k}^{\text{LF}}, \varepsilon_k^{\Diamond, \ell}/2 ), \qquad \forall \ell \in\{0,...,L-1\},
\end{gather*}
from the low-fidelity training set $\Xi_{\Delta, k}^{\text{LF}} = (\bigcup_{\ell=0}^{L-1} \mathcal{E}_k^{P,\ell}) \cup \mathcal{E}_k^{C,L}$.

The absolute tolerances can be derived from the initial ensembles, with a computation which depends on the methodology employed during the analysis step.
In particular, let $\{\varepsilon_r^\ell \in \mathbb{R}^+\}_{\ell = 0}^{L-1}$, $\varepsilon_r^0\leq...\leq \varepsilon_r^{L-1}$, denote the relative tolerances associated with the different accuracy level. The aRB-ML-EnKF absolute tolerances are then computed employing a multi-level estimator of the trace of the state covariance, consistent with Equation \eqref{eq:tolerance}, as
%
\begin{equation*}
\begin{aligned}
\varepsilon^{\text{ML}, \ell}_k = \dfrac{2 \varepsilon_r^\ell}{M_0-1} \sum_{n=1}^{M_0} \| \omega_{k\vert k}^{P, n, 0} - \overline{\omega}_{k\vert k}^{P, 0} \|_\mathcal{V}^2 + \sum_{s=1}^L \dfrac{2\varepsilon_r^\ell}{M_s-1} \sum_{n=1}^{M_s} \Big( \| \omega_{k\vert k}^{P, n, s} - \overline{\omega}_{k\vert k}^{P, s}\|_\mathcal{V}^2 \\
- \| \omega_{k\vert k}^{C, n, s} - \overline{\omega}_{k\vert k}^{C, s}\|_\mathcal{V}^2  \Big) ,
\end{aligned}
\end{equation*}
while the aRB-MF-EnKF tolerances are obtained employing a multi-fidelity estimator, as
\begin{equation*}
\begin{aligned}
\varepsilon_k^{\text{MF}, \ell} = \dfrac{2 \varepsilon_r^\ell / 4^L}{M_0 - 1} \sum_{n=1}^{M_0} \| \omega_{k\vert k}^{P, n, 0} - \overline{\omega}_{k\vert k}^{P, 0} \|_\mathcal{V}^2 + \sum_{s=1}^L \dfrac{2 \varepsilon_r^\ell / 4^{L-s}}{M_s - 1} \sum_{n=1}^{M_s}\Big( \| \omega_{k\vert k}^{P, n, s} - \overline{\omega}_{k\vert k}^{P, s}\|_\mathcal{V}^2 \\ + \dfrac{1}{4} \| \omega_{k\vert k}^{C, n, s} - \overline{\omega}_{k\vert k}^{C, s},\|_\mathcal{V}^2 - ( \omega_{k\vert k}^{P,n,s} - \overline{\omega}_{k\vert k}^{P,s}, \omega_{k\vert k}^{C,n,s} - \overline{\omega}_{k\vert k}^{C,s} )_\mathcal{V} \Big) .
\end{aligned}
\end{equation*}


\section{Numerical experiments on quasi-geostrophic equations}
\label{sec:numerical_experiments}

In this section, we test the proposed techniques on a benchmark problem involving the two-dimensional quasi-geostrophic model. Its non-linearity and high-dimensionality makes it a non-trivial benchmark for data assimilation algorithms. The implementation is done in Python and the simulations were run on an Intel Xeon Platinum 8260. The code is available in the GitLab repository\footnote{\url{https://gitlab.tue.nl/bayanihan/da/enkf/adaptive-hierarchical-enkf}}.

The quasi-geostrophic equations (QGE) describe the motion of stratified, rotating incompressible ideal fluid in the low Rossby number limit and their derivation dates back to the seminal work of Charney \cite{charney1948}. Because of their role in weather prediction and climate modeling, they have been widely investigated in combination with data assimilation and reduced order modeling algorithms, and have been used as benchmark by various authors \cite{evensen1994, popov21, strazzullo2018}.

In this study, we express the QGE directly in the vorticity formulation. For these equations, we introduce two numerical schemes which exploit finite elements or reduced basis methods in space and finite differences in time.
Let $\Omega \subset \mathbb{R}^2$ be a sufficiently regular geometrical domain, $T \in \mathbb{R}^+$ an arbitrary time horizon, and let $\psi, \omega : \Omega \times [0, T] \rightarrow \mathbb{R}$ denote the stream function and the vorticity field, respectively. The stream function is related to the velocity field $\mathbf{u} : \Omega \times [0, T] \rightarrow \mathbb{R}^2$ through the relation
\begin{gather*}
    \mathbf{u} = \left(\dfrac{\partial \psi}{\partial y}, - \dfrac{\partial \psi}{\partial x} \right).
\end{gather*}
Let $\text{Ro} \in \mathbb{R}^+$ be the Rossby number that quantifies the ratio between the inertia and the Coriolis forces; and let $\text{Re} \in \mathbb{R}^+$ be the Reynolds number that represents the ratio between the inertia and the viscous forces. When the Rossby number is sufficiently small, the dynamics of a rotating flow can be described by the quasi-geostrophic equations:
\begin{align}
    \label{eq:strong_QGE_a}
    \dfrac{\partial \omega}{\partial t} - \dfrac{\partial \psi}{\partial x} - \dfrac{\text{Ro}}{\text{Re}} \Delta \omega + \text{Ro} \left( \dfrac{\partial \psi}{\partial y}\dfrac{\partial \omega}{\partial x} - \dfrac{\partial \psi}{\partial x}\dfrac{\partial \omega}{\partial y} \right)  = F, \qquad
    \omega + \Delta \psi = 0.
\end{align}

In our numerical experiments, we consider the domain $\Omega = [0,1]^2$, the forcing term $F=\sin (\pi(y-1))/2$, and the dimensionless quantities $\text{Ro}=0.001$ and $\text{Re}=100$. We close the system of equations \eqref{eq:strong_QGE_a} by imposing zero Dirichlet boundary conditions both on the stream function and on the vorticity field, i.e, we assume $\psi(x,y,t)=\omega(x,y,t)=0$ for all $(x,y)\in\partial\Omega$, $t\in [0,T]$. The vorticity and the stream function initial conditions, $\omega(x,y,0)=\omega_0(x,y)$ and  $\psi(x,y,0)=\psi_0(x,y)$, for all $(x,y)\in\Omega$, are obtained by solving the stationary equation
\begin{align*}
    - \dfrac{\partial \psi_0}{\partial x} - \dfrac{\text{Ro}}{\text{Re}} \Delta \omega_0  = F, \qquad  \omega_0 + \Delta \psi_0 = 0,
\end{align*}
with the same boundary conditions assumed for the non-stationary problem.

The weak formulation of problem \eqref{eq:strong_QGE_a} reads: find $\omega \in C^1([0,T]; L^2(\Omega)) \cap C^0([0,T];$ $ H^1_0(\Omega))$ and $\psi \in C^0([0,T]; H^1_0(\Omega))$ such that for all $\eta, \xi \in H_0^1(\Omega)$
\begin{align}
    \label{eq:weak_QGE_a}
    \int_\Omega \left(\dfrac{\partial \omega}{\partial t} \eta  - \dfrac{\partial \psi}{\partial x} \eta  + \dfrac{\text{Ro}}{\text{Re}} {\nabla} \omega \cdot {\nabla} \eta + \text{Ro} \, J(\omega, \psi, \eta)\right)  \,d\Omega&= \int_\Omega F,\\
    \label{eq:weak_QGE_b}
    \int_\Omega \left(\omega \, \xi - {\nabla} \psi \cdot {\nabla} \xi\right) \,d\Omega &= 0,
\end{align}
where the nonlinear term $J$, defined as
\begin{align*}
    J(\omega, \psi, \eta) =  \psi\dfrac{\partial \omega}{\partial y} \dfrac{\partial \eta}{\partial x} - \psi\dfrac{\partial \omega}{\partial x}\dfrac{\partial \eta}{\partial y},
\end{align*}
is obtained integrating by parts the term $(\mathbf{u}\cdot\nabla\omega)\eta$ and using the incompressibility of the flow.

A `high-fidelity' semi-discretization in space is obtained by restricting the test and trial space of \eqref{eq:weak_QGE_a} and \eqref{eq:weak_QGE_b} to the space of first order Lagrange functions supported on a sufficiently refined triangulation of the geometrical domain. Specifically, in our application, the finite element space, $\mathcal{V} \equiv P^1_h $, has size $\mathcal{N}=3969$ and it is defined over a structured non-uniform grid.

At each iteration, a reduced basis space is constructed employing POD, detailed in Subsection \ref{sec:pod}, on a training set of vorticity fields. This leads to an approximation space $\mathcal{V}_\varepsilon^k := \text{span}\{\phi_1, ..., \phi_{N_\varepsilon^k}\}$ for the vorticity field, and to an approximation space  $\mathcal{U}_\varepsilon^k = \text{span} \{ \zeta_1, ...,  \zeta_{N_\varepsilon^k}\}$ for the stream function. These basis functions of the two spaces satisfy the relationship 
\begin{gather*}
\int_\Omega {\nabla} \zeta_i \cdot {\nabla} \eta \,d\Omega = \int_\Omega \phi_i \, \eta \,d\Omega, \qquad \forall \eta \in P_h^1,\; i \in \{1, ..., N_\varepsilon^k\}.
\end{gather*}

\begin{figure}
    \centering
    \includegraphics[width=0.9\linewidth]{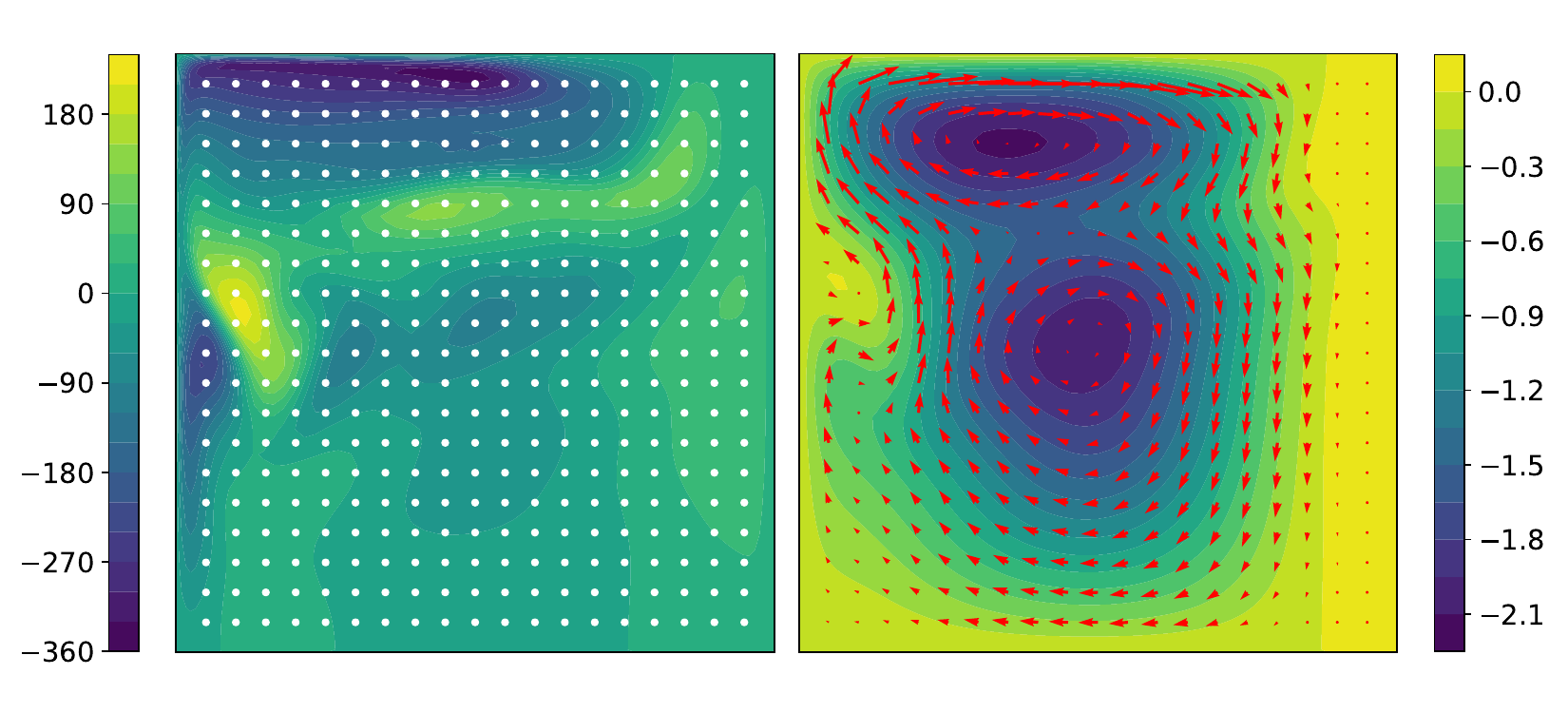}
    \caption{Reference solutions at time $t_0=1000$. On the left: vorticity $\omega(t_0)$ and measurement locations (in white). On the right: stream function $\psi(t_0)$ and corresponding velocity field $\mathbf{u}(t_0)$ (in red).}
    \label{fig:true_solution}
\end{figure}

The semi-discretization gives rise to a set of nonlinear ordinary differential equations that are solved using the implicit midpoint rule. 
The time integration involves dividing the time span $[0,T]$ into subintervals of length $dt = 0.1$. In our context, we set $T=5000$, and we designate the solution over the time range $[1000, 2000]$ as the `true' state trajectory. The solution at time, $t_0=1000$, is shown in Figure \ref{fig:true_solution}.
\firstreview{To characterize the attractor of the dynamics, we consider the system's trajectory over the long term horizon $[4000, 5000]$. The collection of snapshots computed over this interval allows the construction of global reduced basis models. Additionally, randomized selection of these snapshots serves as a sampling strategy, generating states that are distributed according to the invariant measure of the dynamics.}
 
To complete the setup of the data assimilation problem, we introduce the assimilation times $t_k = 1000 + k$, for $k=\{0,...,1000\}$, and the measurement operator $\mathcal{L}:P^1_h \rightarrow \mathbb{R}^{361} \equiv \mathcal{Y}$ that provides pointwise evaluations of the vorticity fields at $361$ spatial locations $\mathbf{x}_{i,j} = (x_i, y_j) \in \Omega$, where $x_i = i/20$ and $y_j = j /20$ for all $i,j \in \{1, ..., 19\}$.  Each noise-free measurement is polluted with i.i.d. Gaussian noise with zero mean and covariance $\sigma \in \mathbb{R}^+$, resulting in the global noise covariance $\boldsymbol{\Sigma} = \sigma^2 \mathbf{I} \in \mathbb{R}^{361 \times 361}$. The measurement process produces a sequence of noise-polluted measurements, $\mathbf{d}_k = \mathcal{L} \omega(t_k) + \boldsymbol{\eta}_k$ with $\boldsymbol{\eta}_k \sim \mathcal{N}(0, \boldsymbol{\Sigma})$, which is used to test different data assimilation algorithms. 

In all the numerical experiments, we employ a noise magnitude of $\sigma = 10^{-4}$ corresponding, for this problem, to a low noise condition. This choice allows for a clear observation of the convergence properties of the data assimilation algorithms over the given assimilation horizon.

\subsection{Numerical results}
\label{sec:numerical_results}

In this section we compare the performances of different data assimilation algorithms on different test settings. Each test involves reconstructing the `true' state trajectory using measurements contaminated by noise. To assess the accuracy of the reconstruction, we monitor the relative error in the reconstruction of the vorticity field. This error is computed by comparing the average vorticity field over the principal ensemble with the `true' state vorticity field, namely
\begin{gather*}
\text{err}(t_k\vert \mathbf{d}_r) = \dfrac{\| \overline{\omega}^P_{k\vert r} - \omega^\star(t_k)\|_{L^2(\Omega)}}{\|\omega^\star(t_k)\|_{L^2(\Omega)}}, \quad r\in\{k, k+1\}, \quad \forall k=\{0,...,1000\}.
\end{gather*}
The error calculation is repeated for four trajectory reconstructions employing independent noise realizations. The resulting confidence intervals enable the distinction between the stochastic fluctuations and the mean behavior of the filters.
As further performance indicators, we monitor how the size of the reduced model changes over the assimilation time and \firstreview{the wall clock time required for the algorithm to execute each data assimilation step, i.e., the time required to complete one prediction and one update step. For the post-processing, we also consider the cumulative computational cost, i.e., the sum of the individual costs associated with all the data assimilation steps up to a given assimilation time, and the long-term computational cost, i.e. the average cost of the last ten data assimilation steps before a given assimilation time.}

As a first experiment, we assess the performance of the standard \firstreview{ensemble Kalman filter} against those of the multi-level and multi-fidelity ensemble Kalman filters \firstreview{under identical long-term computational cost conditions}. For the multi-level and multi-fidelity filters, we consider two levels of accuracy and we employ different approaches to construct surrogate models. Notably, we introduce a reference multi-level filter and a reference multi-fidelity filter that employ the high-fidelity model in place of the surrogate model. Although the resulting filter is prohibitively expensive for practical applications, it does not introduce approximation errors and thus it serves as a reference for the other multi-accuracy algorithms. The considered multi-level and multi-fidelity algorithms encompass the inflation-deflation aRB-ML-EnKF and aRB-MF-EnKF methods, and a filter that makes use of memory-less reduced order modeling, detailed in \cite{donoghue22} and already mentioned in Remark \ref{rmk:yano}.


All filters employ ensembles initialized with states sampled from a Gaussian distribution with inverse Laplacian covariance, as detailed in Subsection \ref{sec:initial_ensemble}. \firstreview{The standard EnKF ensemble has $N_P = 256$ high-fidelity states}, whereas the multi-accuracy filters have $N_P = 48$ principal states, $N_C = 48$ control states, and $N_A = 1200$ ancillary states. The relative accuracy for the surrogate models is set to $\varepsilon_r = 10^{-3}$.
\begin{figure}
    \includegraphics[width=\linewidth]{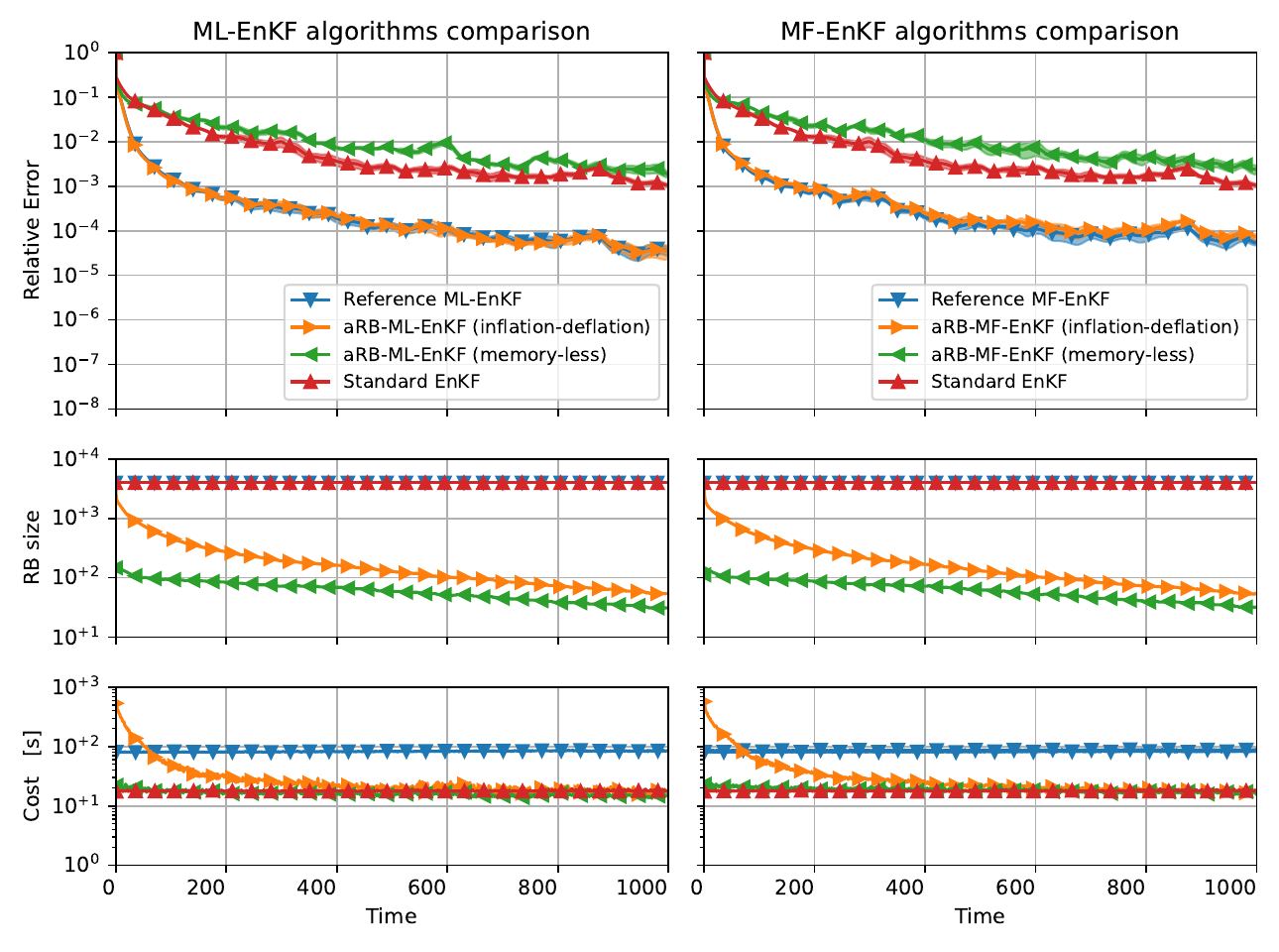}
    \caption{\firstreview{This plot compares different data assimilation algorithms. The left column shows the results of filters employing multi-level updates, while  the right column shows those using multi-fidelity updates.
    The top row reports the relative errors versus assimilation time, the center row displays how the reduced basis size changes over the assimilation time, and the bottom row shows the wall clock time per single data assimilation step versus assimilation time. }
    }
    \label{fig:comparison_algorithms}
\end{figure}

\begin{figure}
    \includegraphics[width=\linewidth]{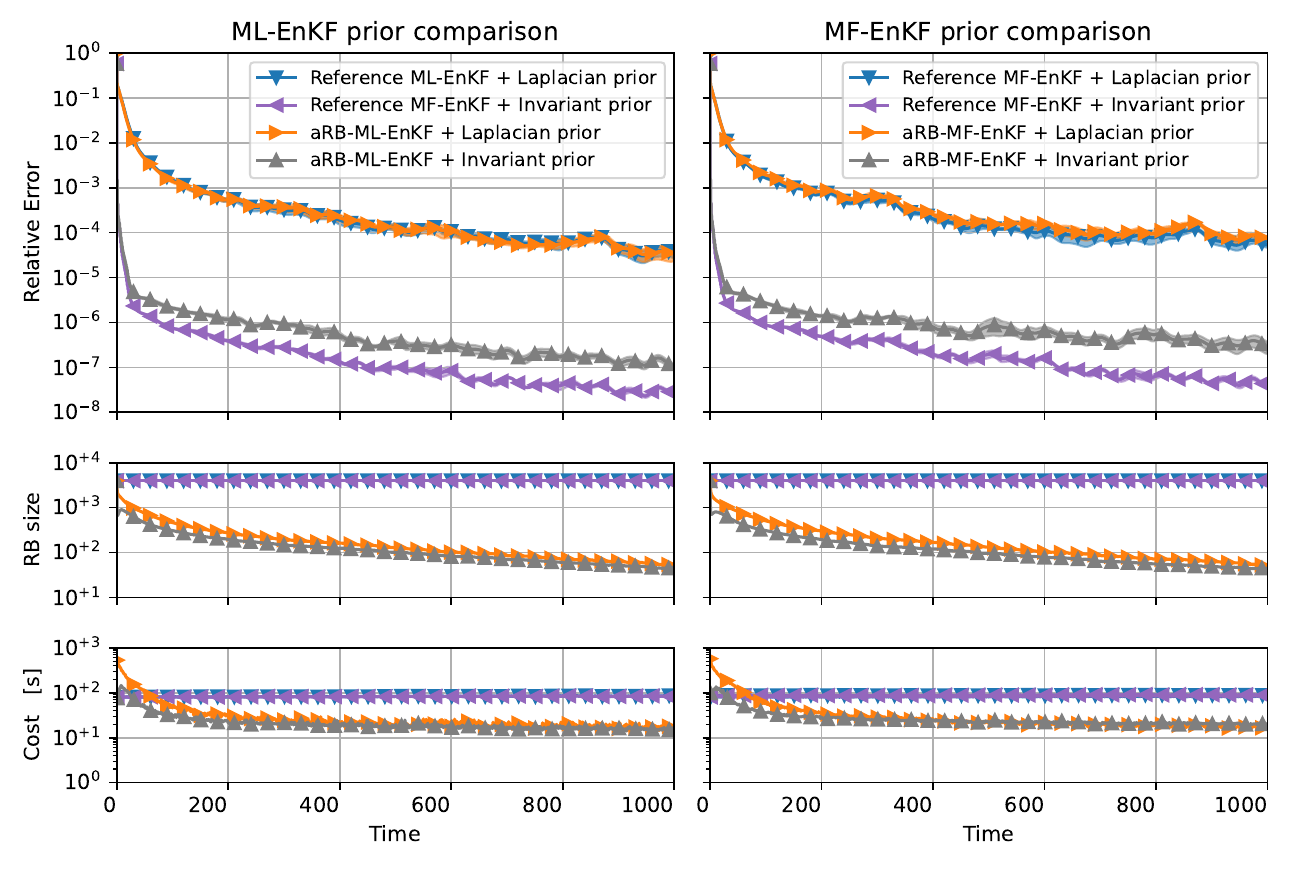}  
    \caption{This plot compares different sampling priors. The left column shows the results of filters employing multi-level updates, while  the right column shows those using multi-fidelity updates. The sub-plots in the top row illustrate relative errors versus assimilation time, those in the center row display the relationship between reduced basis size and assimilation time, and those at bottom the \firstreview{wall clock time per single data assimilation step} versus assimilation time. 
   }
    \label{fig:comparison_priors}
\end{figure}

In Figure \ref{fig:comparison_algorithms}, we present the results of this experiment. The error plots show good performances of the inflation-deflation aRB-ML-EnKF and aRB-MF-EnKF, whose error closely aligns to the error of the reference algorithm, both in the multi-level and in the multi-fidelity case, and \firstreview{it is 1 to 1.5 orders of magnitude smaller} than the errors obtained with the standard EnKF and with the memory-less algorithm.
The evolution of the size of the surrogate model reveals that the inflation-deflation approach ensures larger and more accurate reduced spaces compared to the memory-less algorithm. 
\firstreview{However, this comes with a higher computational cost in the initial assimilation windows, which can be mitigated by selecting a better prior distribution, as shown in Figure \ref{fig:comparison_priors}.}
Finally, we remark that, for the same computational cost, the multi-level update produces slightly more accurate reconstructions than the  multi-fidelity one, although the behavior of the two algorithms does not differ significantly.

In a second experiment, we compare the reconstructions obtained employing a Gaussian prior with inverse Laplacian covariance, which encodes $H^1(\Omega)$ regularity, and the reconstructions achieved using the invariant measure of the dynamics as a prior distribution. Here we focus our analysis on two algorithms: the inflation-deflation aRB-ML-EnKF/aRB-MF-EnKF method, and the reference algorithm employing the high-fidelity model. The ensemble sizes and relative accuracy are taken as in the previous experiment: $N_P = 48$, $N_A = 1200$, and $\varepsilon_r = 10^{-3}$.

The results of the experiment, summarized in Figure \ref{fig:comparison_priors}, demonstrate the benefits of integrating additional information about the system dynamics into the data assimilation process. Using the invariant measure of dynamics as a prior leads to significantly faster convergence for both algorithms. However, it also poses a challenge for the adaptive reduced order modeling scheme. Specifically, the error plot reveals over the last assimilation windows a one order of magnitude difference in relative errors between the inflation-deflation aRB-ML-EnKF/aRB-MF-EnKF method, and the reference algorithm. 
To better understand this discrepancy, we can examine the plots depicting the surrogate model size (central row of Figure \ref{fig:comparison_priors}). There, we observe a rapid contraction of the approximation space during the first \firstreview{data assimilation steps}, suggesting a degradation in the model's approximation properties. This information loss is mitigated by the inflation-deflation approach, which extracts additional information from the system dynamics, thereby expanding the size of the approximation space. Nevertheless, this compensatory mechanism does not appear sufficient to close the gap that arises in the early \firstreview{data assimilation steps}. 

In the next experiment, we compare the state reconstructions obtained with surrogate models with different relative error tolerances. We consider the reference filter and the inflation-deflation aRB-MF-EnKF/aRB-ML-EnKF, with three possible values for the relative tolerance, $\varepsilon_r \in \{10^{-3}, 10^{-2}, 10^{-1}\}$. We employ initial ensembles sampled from a Gaussian prior with inverse Laplacian covariance, and we compare results obtained with multi-level and with multi-fidelity update steps. The number of states in the ensembles is the same used in the previous experiments: $N_P = 48$,  $N_C = 48$, and $N_A=1200$.
\begin{figure}
    \includegraphics[width=\linewidth]{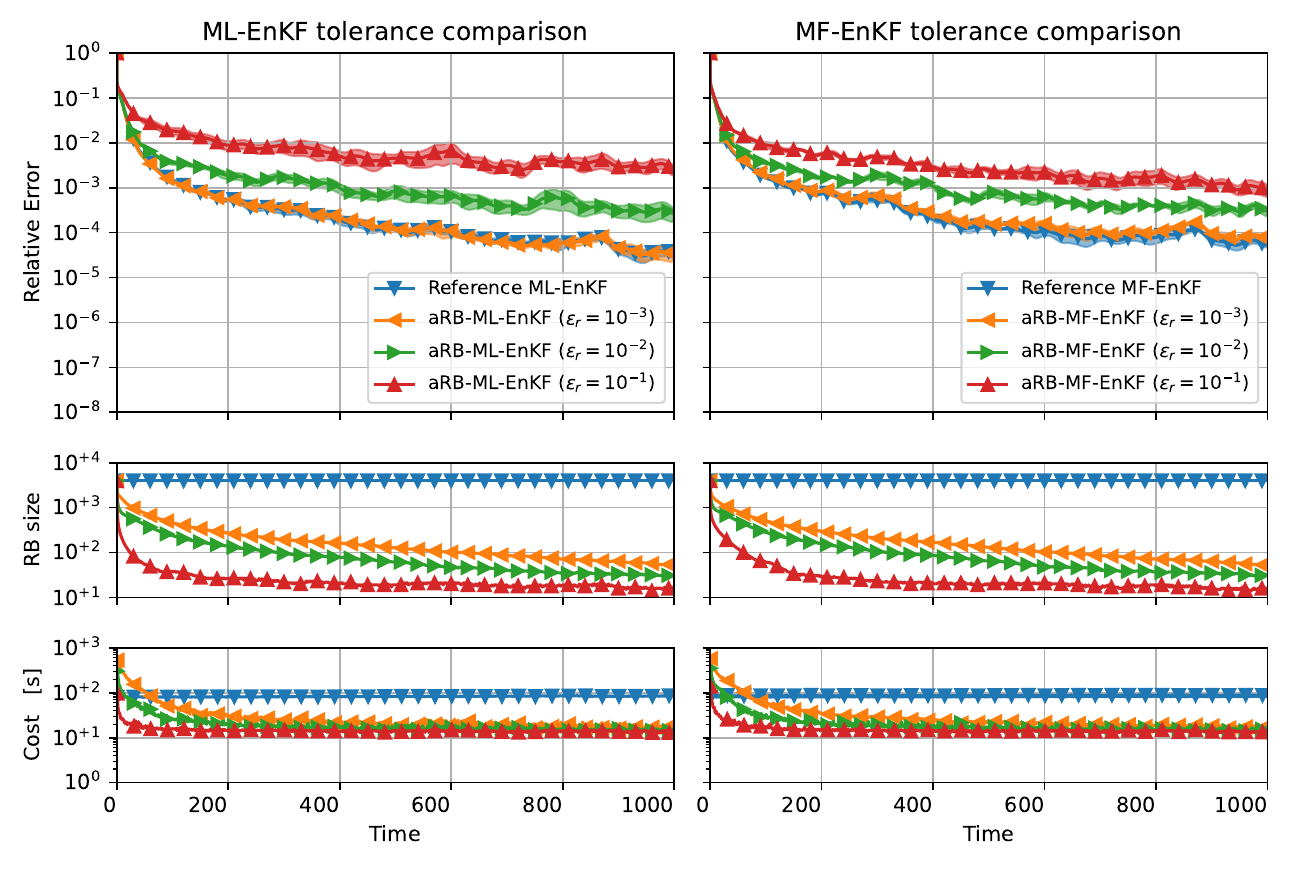}    
    \caption{This plot compares different tolerances for the surrogate model. The left column shows the results of filters employing multi-level updates, while  the right column shows those using multi-fidelity updates. The sub-plots in the top row illustrate relative errors versus assimilation time, those in the center row display the relationship between reduced basis size and assimilation time, and those at bottom the \firstreview{wall clock time per single data assimilation step} versus assimilation time. 
   }
    \label{fig:comparison_tolerance}
\end{figure}

The results, depicted in Figure \ref{fig:comparison_tolerance}, show a rapid convergence in the reconstruction error for both the multi-level and multi-fidelity algorithms as the model tolerance decreases. Already with a relative tolerance of $\varepsilon_r = 10^{-3}$, the introduced approximation error is negligible compared to other sources of error. Notably, the sensitivity to the approximation error differs between the two approaches: the multi-fidelity approach exhibits lower sensitivity, resulting in lower reconstruction errors when larger tolerances are considered. Conversely, it is the multi-level approach that achieves the lowest error.
%

\begin{figure}
    \includegraphics[width=\linewidth]{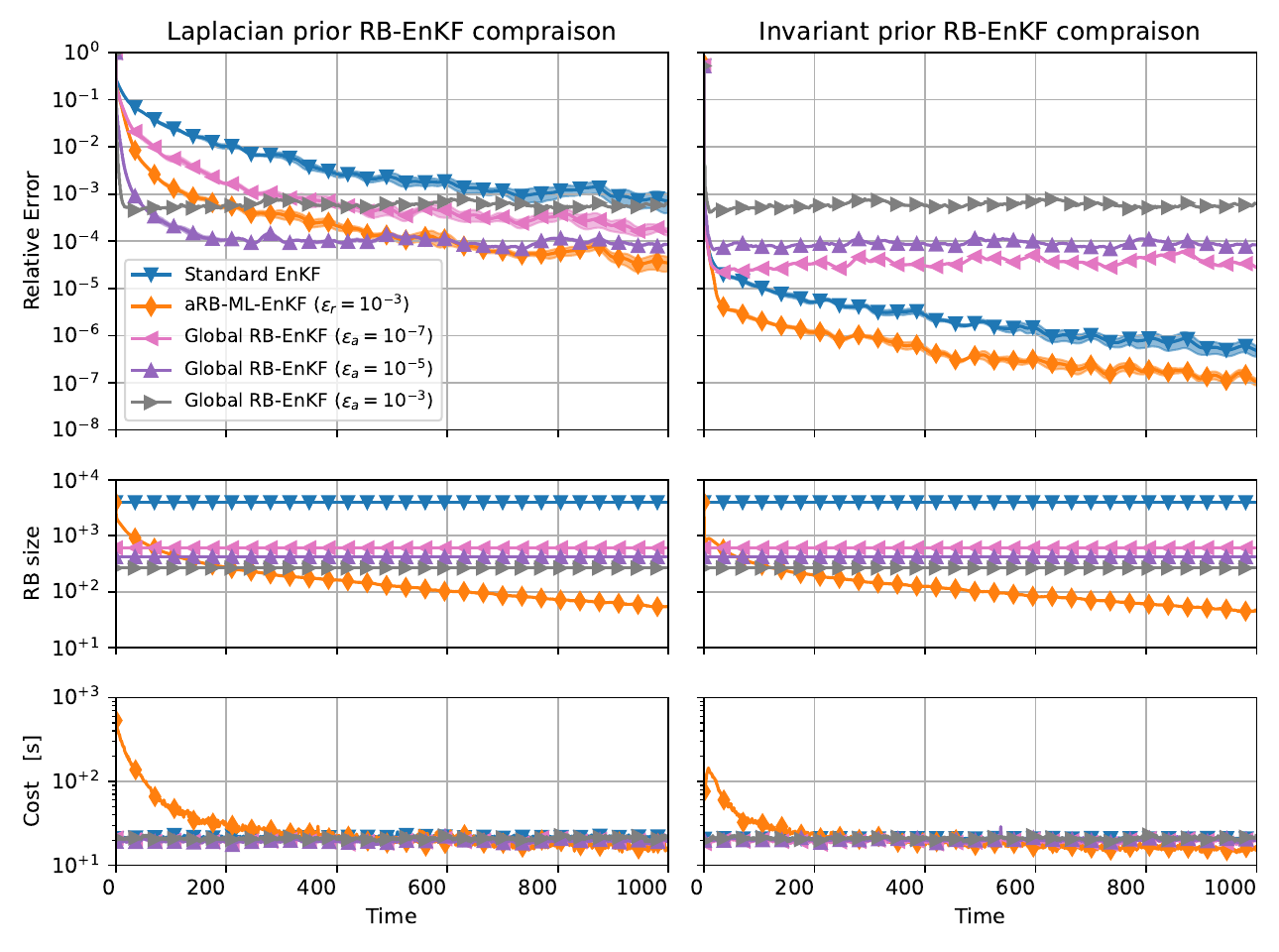}    
    \caption{\firstreview{This plot compares different data assimilation algorithms: the standard EnKF, the two-level aRB-ML-EnKF, and three RB-EnKF with different accuracy levels and ensemble sizes. In the left column are reported the results of filters employing the inverse Laplacian prior, while in the right column those using the invariant measure prior. The top row reports the relative errors versus assimilation time, the center row displays how the reduced basis size changes over the assimilation time, and the bottom row shows the wall clock time per single data assimilation step versus assimilation time.}}
    \label{fig:comparison_rb-enkf}
\end{figure}

In a fourth experiment, we compare the state reconstruction obtained with the standard EnKF and with the aRB-ML-EnKF to those obtained using a RB-EnKF that employs global reduced basis models of different accuracy. These models are constructed offline, leveraging the same set of full order solutions used to obtain samples from the invariant measure of the dynamics. For the aRB-ML-EnKF, we consider two levels of accuracy, while for the other filters, we use a single level of accuracy. We use initial ensembles sampled from a Gaussian distribution with inverse Laplacian covariance and from the same invariant measure used for the offline construction of the RB-EnKF surrogate models. The size of all ensembles is chosen to ensure the same long-term computational cost per data assimilation step. Specifically, for the standard EnKF we consider ensembles of $N_P = 256$ high-fidelity states, whereas for the multi-level filter, we use a principal ensemble of size $N_P = 48$, an ancillary ensemble of size $N_A = 1200$, and set the relative accuracy to $\varepsilon_r = 10^{-3}$. For the global RB-EnKF, we employ three surrogate models with absolute accuracy $\varepsilon_a \in \{ 10^{-7},10^{-5},10^{-3} \}$ corresponding to low-fidelity ensembles of size $N_P \in \{400, 800, 1600\}$, respectively.%

In Figure \ref{fig:comparison_rb-enkf},  we present the results of the experiment. Using low-fidelity models in place of high-fidelity models significantly reduces the convergence time of the EnKF, especially when the surrogate models have lower accuracy. This effect can be explained by the fact that the global reduced basis models retain information about the system dynamics, allowing the filter to operate within a lower dimensional approximation space. However, for both priors used, this faster convergence has a limit on achievable accuracy, with less accurate surrogate models leading to higher reconstruction errors. This error stagnation aligns with the accuracy of the approximation spaces considered. When comparing these RB-EnKF filters to our proposed adaptive approach, the adaptive method appears to be more accurate in the long run. Specifically, in the case of an invariant measure prior, where all methods benefit from knowledge of the system's long-term behavior, both the standard EnKF and the aRB-ML-EnKF outperform the three global RB-EnKF methods in terms of achieved reconstruction error, with the aRB-ML-EnKF being the most accurate.

\firstreview{
\begin{figure}
    \includegraphics[width=\linewidth]{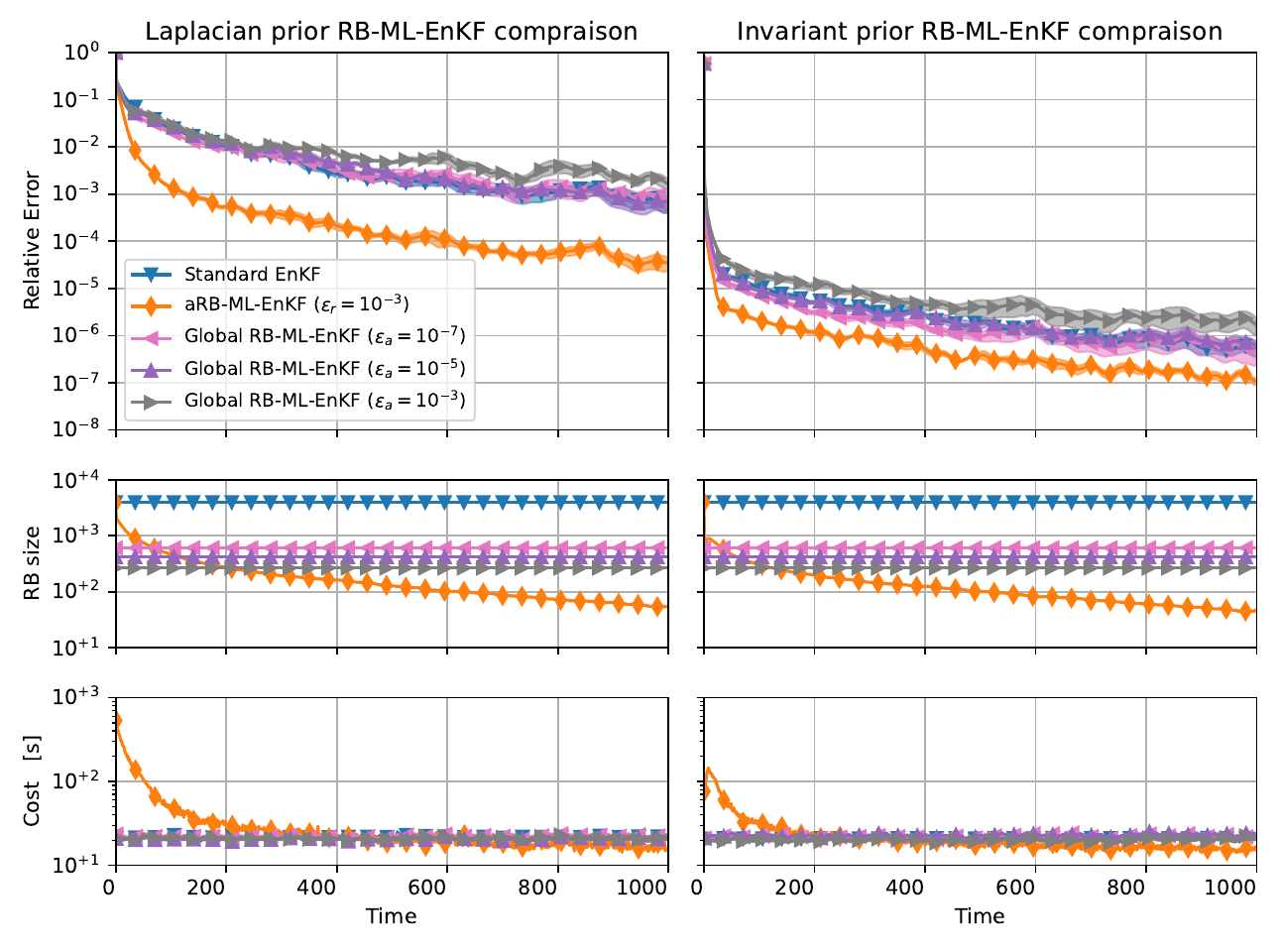}    
    \caption{\firstreview{This plot compares different data assimilation algorithms: the standard EnKF, the two-level aRB-ML-EnKF, and three two-level RB-ML-EnKF with offline-trained low-fidelity models. The left column shows the results of filters employing the inverse Laplacian prior, while in the right column those using the invariant measure prior. The top row reports the relative errors versus assimilation time, the center row displays how the reduced basis size changes over the assimilation time, and the bottom row shows the wall clock time per single data assimilation step versus assimilation time.}}
    \label{fig:comparison_rb_ml_enkf}
\end{figure}

In the next experiment, we compare the state reconstruction obtained with the standard EnKF and with the aRB-ML-EnKF to those obtained from three multi-level ensemble Kalman filters employing the same global low-fidelity models used in the previous experiment. We consider two levels of accuracy for all the filters except the standard EnKF, and we use the same inverse Laplacian and measure invariant priors used before. The size of the ensembles is again chosen to ensure the same long-term computational cost for the filters. For the standard EnKF and for the aRB-ML-EnKF, we consider the same testing conditions as from the previous experiment, while for the global RB-ML-EnKF, we employ principal ensembles of full order states of size $N_P=48$, and ancillary ensembles of size $N_A \in \{ 330, 660, 1320 \}$, respectively associated to the global reduced basis model with accuracy $\varepsilon_a \in \{ 10^{-7}, 10^{-5}, 10^{-3} \}$.

The results of this experiment are presented in Figure \ref{fig:comparison_rb_ml_enkf}. Contrary to the previous test, the convergence of the global RB-ML-EnKF algorithms does not seem to be affected by the accuracy of the global low-fidelity models. However, the accuracy achieved with the non-adaptive multi-level filters is comparable to that of the standard EnKF at the same cost, i.e., one to two orders of magnitude less accurate than the reconstruction obtained with the proposed adaptive method. 
}

\firstreview{In the sixth experiment}, we compare the state reconstructions obtained employing the inflation-deflation ML-EnKF with ensembles of varying sizes. Specifically, we explore the impact of the principal and of the ancillary ensemble size on the reconstruction error, on the surrogate model size, and on the computational cost. 
For this analysis, we first consider a scenario where the principal and control ensemble sizes are fixed, $N_P = N_C = 48$, while the ancillary ensemble size $N_A$ varies in $\{600, 1200, 2400, 4800\}$. Subsequently, we explore the scenario where the ancillary ensemble size is fixed, $N_A = 1200$, and the principal and control ensemble sizes vary, $N_P=N_C \in \{12, 48, 192, 768\}$.
In both cases, we maintain a relative tolerance of $\varepsilon_r = 10^{-3}$ and we employ initial ensembles sampled from a Gaussian prior with inverse Laplacian covariance.
\begin{figure}
    \includegraphics[width=\linewidth]{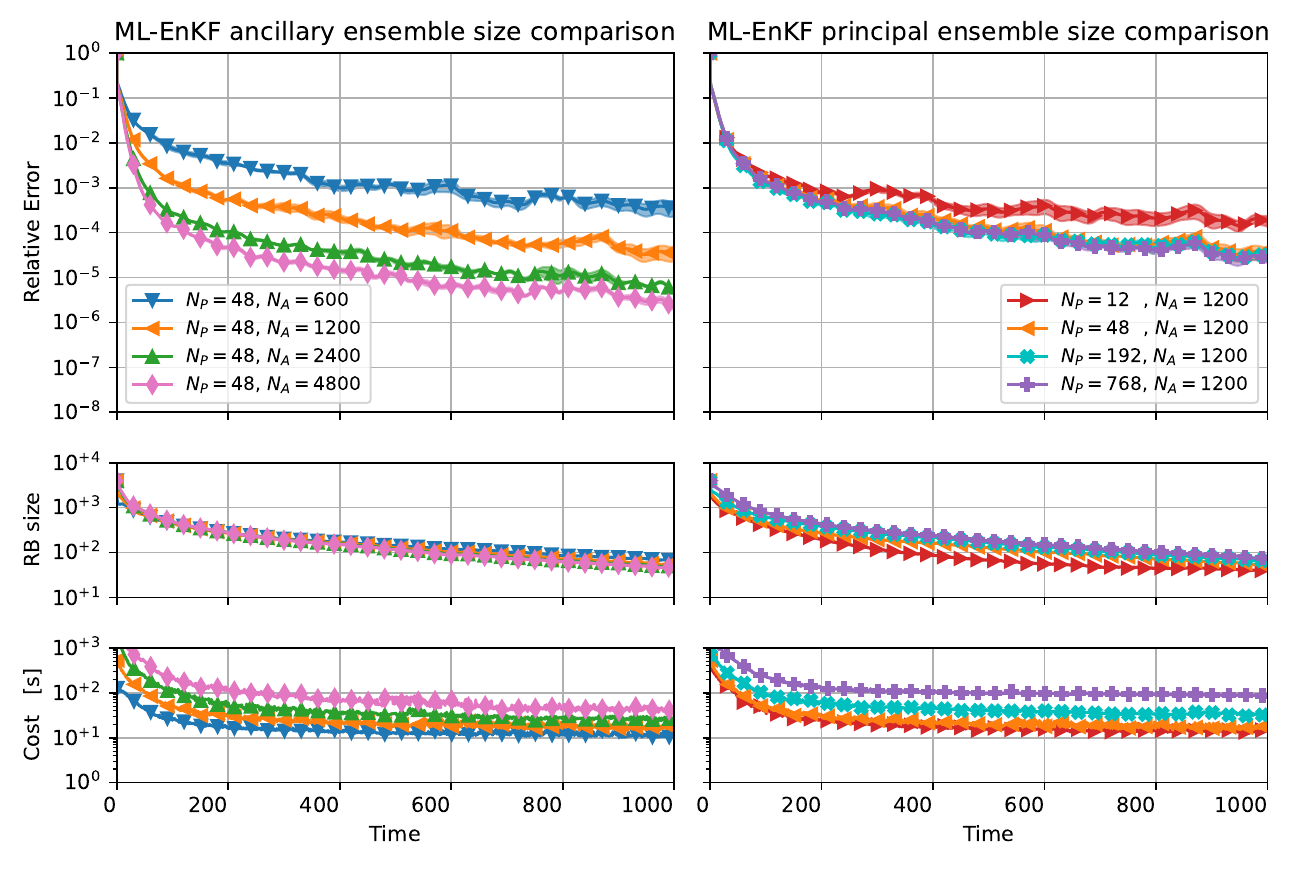}    
    \caption{This plot compares the results of the inflation-deflation ML-EnKF algorithm obtained by employing ensembles of different sizes. In the left column, varying values of the ancillary ensemble size are considered, while in the right column, varying values of the principal and control size are considered. The sub-plots in the top row illustrate relative errors versus assimilation time, those in the center row display the relationship between reduced basis size and assimilation time, and those at bottom the \firstreview{wall clock time per single data assimilation step} versus assimilation time. 
    }
    \label{fig:comparison_size}
\end{figure}

The results of the experiment with fixed principal and control ensemble size are summarized in Figure \ref{fig:comparison_size}, on the left, and show a reduction in the reconstruction error with the increase of the ancillary ensemble size. This corresponds to a higher computational cost and to an almost constant trend for the surrogate model size. This result implies that despite the surrogate model effectively approximates the full order model, the filter is still distant from its mean-field behavior. In such circumstances, a larger ancillary ensemble or a larger number of levels could be used to improve the state reconstruction.

In contrast, the experiment involving a fixed ancillary ensemble size, illustrated in Figure \ref{fig:comparison_size} on the right, shows an increase in the reconstruction error for smaller sizes of the principal and control ensembles. This increase, associated with a reduction in the surrogate model size, is a consequence of insufficient information provided by the principal ensemble during the inflation phase. The same behavior is not observed for larger ensemble sizes, where the expansion of the principal ensemble only increases the computational cost of the method without yielding any improvement in the reconstruction error. We must mention that although a substantial expansion of the principal ensemble could shift the filter towards its mean field limit and thus improve the reconstruction error, such an expansion would come at a prohibitively high cost.

\begin{figure}
    \includegraphics[width=\linewidth]{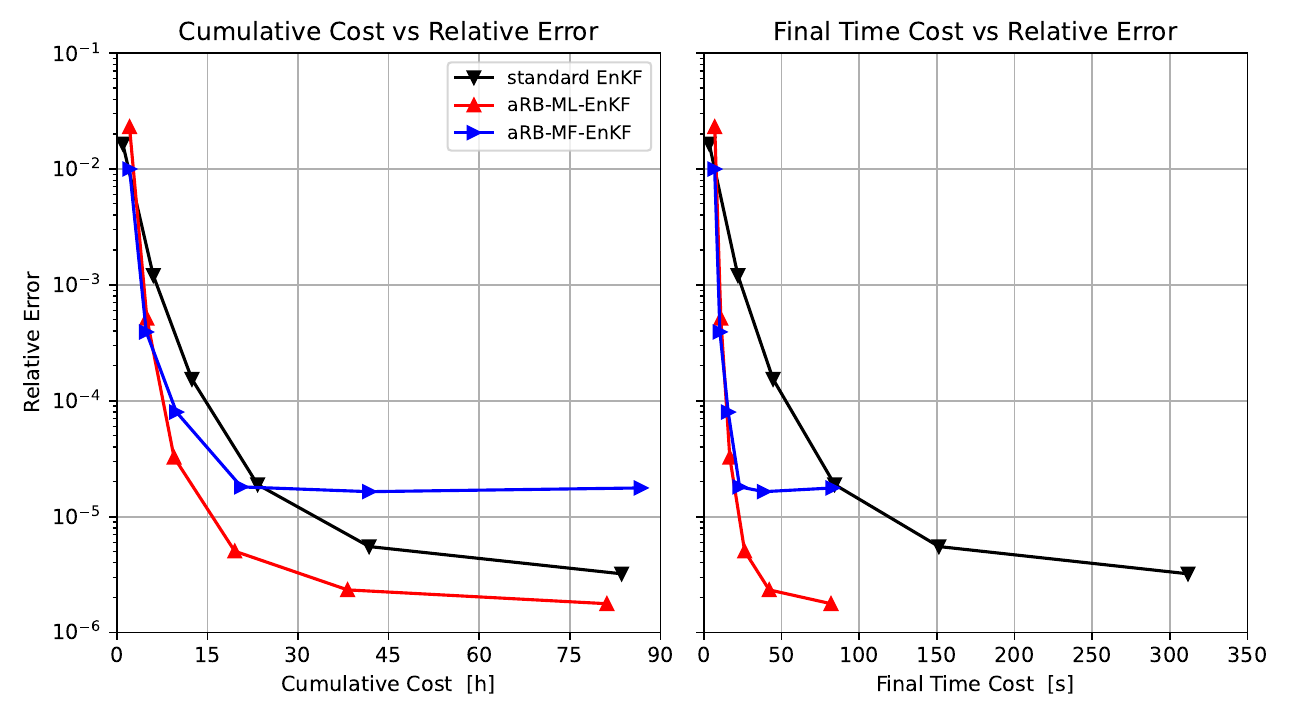}    
    \caption{\firstreview{The plot compares the reconstruction error and computational cost of three different data assimilation algorithms. On the left: the relative error versus the cumulative clock time required for all data assimilation steps. On the right: the relative error versus the clock time of the last data assimilation step.}}
    \label{fig:comparison_cost}
\end{figure}

In the last experiment, we compare the accuracy and computational cost of the standard EnKF with those of the inflation-deflation aRB-ML-EnKF and aRB-MF-EnKF. For the standard EnKF, we consider ensembles of size $N_P \in \{ 48, 300, 600, 1200, 2400, 4800 \}$, while for the multi-level and multi-fidelity filters, we use two levels of accuracy, with a principal ensemble of fixed size $N_P = 48$, and ancillary ensembles of size $N_A \in \{ 48, 600, 1200, 2400, 4800, 9600 \}$. The states in the ensemble of the standard EnKF and in the principal ensemble of the multi-fidelity and multi-level EnKF are forecasted using the high-fidelity model.
The results, presented in Figure \ref{fig:comparison_cost}, display the reconstruction error after $1000$ data assimilation steps, compared to the cumulative cost of the algorithms up to the final time $t_{1000}$ (left plot), and the cost of a single data assimilation step at the final time (right plot). The aRB-ML-EnKF appears to be the least costly of the methods, approaching the mean-field limit faster than the standard EnKF. For example, to achieve an accuracy of $10^{-5}$ at the final time, the standard EnKF requires a cumulative computational cost that is double that of the aRB-ML-EnKF.
In contrast, the aRB-MF-EnKF exhibits a more complex behavior: initially, with an ancillary ensemble of less than 2400 states, the algorithm is more accurate than the standard EnKF for the same cost. However, for larger ensemble sizes and higher cumulative costs, its performance significantly deteriorates. This behavior, rather than indicating slower convergence to the same mean-field limit, seems to align with the predictions of Remark \ref{remark_1}: the MF-EnKF does not appear to share the same mean-field limit as the other two filters, making it, at least in this example, a less effective filter.

\section{Conclusions}
\label{sec:conclusions}

In conclusion, we propose two methods to effectively incorporate reduced-order modeling techniques within the multi-level and multi-fidelity extensions of ensemble Kalman filtering. The on-the-fly construction of reduced basis models, using the suggested inflation-deflation approach, proves highly effective in generating low-dimensional surrogate models that are both inexpensive to assemble and assess, while maintaining good approximation properties. This is particularly relevant for problems that are globally irreducible due to initial uncertainty and the complexity of their dynamics. 

\firstreview{
The numerical experiments show that the adaptive multi-level algorithm achieves errors comparable to those obtained using the full order model alone, yet at a significantly lower cost.
Furthermore, since the reconstruction error does not deviate from the reference behavior as the state uncertainty decreases, we can conclude that the proposed adaptive tolerance is effective in ensuring that surrogate models maintain the necessary accuracy.
Between the aRB-ML-EnKF and the aRB-MF-EnKF, the former appears to be a more reliable option. Indeed, although in most of the experiments the adaptive multi-fidelity approach yields similar results, the last experiment involving larger ensembles seems to suggest convergence to a different, and in this case less accurate, mean field limit. 
}

Regarding future developments, integrating adaptive hyper-reduction techniques into the on-the-fly training phase of the surrogate model, as in \cite{donoghue22}, is expected to enhance the algorithm's performance and handle non-linearity more efficiently. This, together with the implementation of the method using more than two levels of accuracy, constitutes the natural progression of the proposed approaches.

\bibliography{references}

\end{document}